\documentclass[a4paper, 12pt]{extarticle}
\usepackage{geometry}
\usepackage{graphics}
\usepackage{dcolumn}
\usepackage{bm}
\usepackage{amsmath}
\usepackage{hyperref}
\hypersetup{colorlinks=true, urlcolor=blue}
\usepackage{color}
\usepackage{tikz}
\usepackage{multirow}
\usepackage{cite}
\usepackage{float}
\usepackage[caption = false]{subfig}
\usepackage{newlfont}
\usepackage{amssymb}
\usepackage{amsfonts}
\usepackage{amsthm}
\usepackage{graphicx}
\usepackage{graphics,wrapfig,times}
\usepackage{bm}

\def \be{\begin{equation}}
\def \ee{ \end{equation} }
\def \bea{\begin{eqnarray}}
\def \eea{\end{eqnarray}}

\begin{document}

\definecolor{red}{rgb}{1,0,0}
\title{Additive Sequences, Sums, Golden Ratios\\ and Determinantal Identities}
\author{Asutosh Kumar\\
\small P. G. Department of Physics, Gaya College, Magadh University, Rampur, Gaya 823001, India\\ \small Harish-Chandra Research Institute, HBNI, Chhatnag Road, Jhunsi, Prayagraj 211 019, India\\ \small Vaidic and Modern Physics Research Centre, Bhagal Bhim, Bhinmal, Jalore 343029, India\\ \small (asutoshk.phys@gmail.com)}
\date{}

\maketitle

\textbf{Abstract.} 
The Fibonacci sequence is a series of positive integers 
in which, starting from \(0\) and \(1\), every number is the sum of two previous numbers, and the limiting ratio of any two consecutive numbers of this sequence is called the {\it golden ratio}. The Fibonacci numbers and the golden ratio are two significant concepts that keep appearing everywhere. In this article, we investigate the following issues:
\begin{itemize}
\item We recall the Fibonacci sequence, the golden ratio, their properties and applications, and some early generalizations of the golden ratio. 
The Fibonacci sequence is a \(2\)-sequence because it is generated by the sum of two previous terms, $f_{n+2} = f_{n+1} + f_{n}$. As a natural extension of this, we introduce several typical \(p\)-sequences where every term is the sum of {\it \(p\) previous} terms given \(p\) initial values called {\it seeds}. In particular, we introduce the notion of \(1\)-sequence. We then discuss generating functions and limiting ratio values of \(p\)-sequences.
Furthermore, inspired by Fibonacci's rabbit pair problem, we consider a general problem whose particular cases lead to nontrivial additive sequences.

\item We obtain closed expressions for odd and even sums, sum of the first \(n\) numbers, and the sum of squares of the first \(n\) numbers of the {\it exponent} \(p\)-sequence whose {\it seeds} are $(0,1,\cdots,p-1)$. 

\item We investigate the \(p\)-golden ratio of \(p\)-sequences,  express a positive integer power of the \(p\)-golden ratio as a polynomial of degree $p-1$, and obtain values of golden angles for different \(p\)-golden ratios. We also consider further generalizations of the golden ratio.  
%

\item We establish a family of determinantal identities of which the Cassini's identity is a particular case.
\end{itemize}

\section{Introduction}

We are familiar with the celebrated Fibonacci sequence \cite{vorobyov1963, hoggatt1969, koshy2001}: 1, 1, 2, 3, 5, 8, 13, 21, 34, 55, ... . In this sequence, each number is the sum of the previous two, starting from 1 and 1. The ratio of consecutive Fibonacci numbers approaches the unique number 1.618. That is, $f_n = f_{n-1} + f_{n-2}$ with $f_1 = 1$ and $f_2 = 1$
\footnote{One can also start with $f_0 = 0$ and $f_1 = 1$. In that case the Fibonacci sequence will be $\{0, 1, 1, 2, 3, 5, 8, 13, 21, 34, 55, ... \}$.},
and $\lim_{n \rightarrow \infty} \frac{f_{n+1}}{f_n} = 1.618$ (upto three decimal places). This sequence arose from the Fibonacci's famous rabbit pair problem. A version of this problem is: {\it A man puts a male-female pair of adult rabbits in a field. Rabbits take a month to mature before mating. One month after mating, females give birth to one male-female pair and then mate again. It is assumed that no rabbits die but continue breeding. How many rabbit pairs are there after one year?} See Table \ref{table-fibonacci} for the answer.

\begin{table}[]
\centering
\begin{tabular}{|l|l|l|l|}
\hline
$n$ (month)  & adult pair $(a_n)$  & baby pair $(b_n)$ & total pair $(t_n)$ \\ \hline
0  & 1        & 0        & 1     \\ \hline
\hline
1  & 1        & 0        & 1     \\ \hline
2  & 1+0      & 1        & 2     \\ \hline
3  & 1+1      & 1        & 3     \\ \hline
4  & 2+1      & 2        & 5     \\ \hline
5  & 3+2      & 3        & 8     \\ \hline
6  & 5+3      & 5        & 13    \\ \hline
7  & 8+5      & 8        & 21    \\ \hline
8  & 13+8     & 13       & 34    \\ \hline
9  & 21+13    & 21       & 55    \\ \hline
10 & 34+21    & 34       & 89    \\ \hline
11 & 55+34    & 55       & 144   \\ \hline
12 & 89+55    & 89       & 233   \\ \hline
13 & 144+89   & 144      & 377   \\ \hline
14 & 233+144  & 233      & 610   \\ \hline
15 & 377+233  & 377      & 987   \\ \hline
16 & 610+377  & 610      & 1597  \\ \hline
17 & 987+610  & 987      & 2584  \\ \hline
18 & 1597+987 & 1597     & 4181  \\ \hline
\end{tabular}%
\caption{The Fibonacci sequence resulting from Fibonacci's famous rabbit pair problem. 
In the beginning $(n = 0)$, there is one male-female adult pair. At the start of the first month, 
there is one adult pair (they mate) and zero juvenile pair so there is only \(1\) rabbit pair.
At the start of the second month they produce a new pair, so there are \(2\) pairs in the field.
At the start of the third month, the original pair produce a new pair, but the second pair only mate without breeding, so there are \(3\) pairs in all.
And so on.
Note that $a_{n \ge 2} = a_{n-1} + b_{n-1}$, $b_{n \ge 2} = a_{n-1}$, and $t_{n \ge 2} = a_{n} + a_{n-1} = t_{n-1} + t_{n-2}$.
}
\label{table-fibonacci}
\end{table}

Another problem, a modified version of Pingala's (c. 200 BC) \cite{weber1863, singh1936, singh1985, nooten1993, wilson1993, wilson2001}, which yields the Fibonacci sequence is: {\it Suppose $\{s_k \equiv k \}_{k=1}^p$ is the set of syllable elements, and a room of \(n\) syllables is available. In how many ways this \(n\)-syllable room can be occupied by these syllable elements?} See Table \ref{table-2syllable} for the answer. Indeed, there are many ways to obtain the Fibonacci sequence.

\begin{table}[]
\centering
\resizebox{\textwidth}{!}{%
\begin{tabular}{|l|l|l|}
\hline
$n$ & possible arrangements                        & total \\ \hline
1   & $s_1^{1}s_2^{0}(1)~[1]$                                            & {\bf 1}  \\ \hline
2   & $s_1^{2}s_2^{0}(1)~[11],~s_1^{0}s_2^{1}(1)~[2]$                                & {\bf 2}  \\ \hline
\hline
3   & $s_1^{3}s_2^{0}(1)~[111],~s_1^{1}s_2^{1}(2)~[12,~21]$                                 & 3  \\ \hline
4   & $s_1^{4}s_2^{0}(1)~[1111],~s_1^{2}s_2^{1}(3)~[112,~121,~211],~s_1^{0}s_2^{2}(1)~[22]$       & 5  \\ \hline
5   & $s_1^{5}s_2^{0}(1)~[11111],~s_1^{3}s_2^{1}(4)~[1112,~1121,~1211,~2111],~s_1^{1}s_2^{2}(3)~[122,~212,~221]$        & 8  \\ \hline
\end{tabular}%
}
\caption{Possible arrangements of occupation of \(n\)-syllable room with 1-syllable and 2-syllable elements. Each occupation is of the form $s_1^{n_1}s_2^{n_2}(m)$, where $s_k^{n_k}$ means that $s_k$ occurs $n_k$ times and number \(m\) in the parenthesis denotes the possible arrangements or {\it multiplicity} of $s_1^{n_1}s_2^{n_2}$. The numbers in the last column build up a sequence. We call this {\it syllable} \(2\)-sequence.
This can be straightforward generalized for any number of syllable elements: $\{s_k \equiv k \}_{k=1}^p$. Note that $n = \sum_{k=1}^p k~n_k$ and multiplicity $m = \frac{(n_1 + n_2 + \cdots + n_p)!}{n_1! ~n_2! \cdots n_p!}$. And the last column (total number of ways in which \(n\)-syllable room can be occupied by these syllable elements): numbers in the first \(p\) rows will be $2^0,~2^1,~\cdots,2^{p-1}$, and number in the $n^{th}$ row $(n > p)$ will be the sum of \(p\)-previous terms. Numbers $2^0,~2^1,~\cdots,2^{p-1}$ serve as the {\it seeds} for the syllable \(p\)-sequence.}
\label{table-2syllable}
\end{table}

\subsection{Historical background}
\begin{figure}%
\centering
\includegraphics[width = 4.0in]{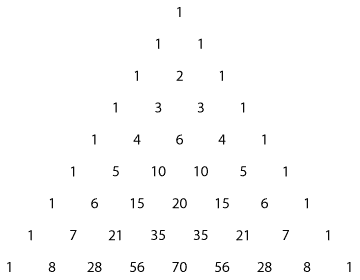}%
\caption{Pingala's Mount Meru (Pascal's triangle).}
\label{fig-pingala}
\end{figure}

It is acknowledged that the notions of binomial coefficients via the {\it Mount Meru} and the Fibonacci sequence were well known to Indian mathematicians--Pingala (c. 200 BC), Varahamihira (505-587), Kedara ($7^{th}$ century), Virahanka ($7^{th}$ century), Halayudha ($10^{th}$ century), Gopala (c. 1135) and Hemachandra (1089-1172) \cite{weber1863, singh1936, singh1985, nooten1993, wilson1993, wilson2001}, and Persian mathematicians--Al-Karaji (953-1029) and Omar Hayyam (1048-1131) (see \cite{coolidge1949}) before Fibonacci who had introduced it to the Western world in his book {\it Liber Abaci} (1202)
\footnote{We, therefore, advocate that the Fibonacci sequence be called the Pingala sequence or the Pingala-Fibonacci sequence.}. 
The shallow diagonals of the Mount Meru sum to the Fibonacci numbers (see Fig. \ref{fig-pingala}), and the Mount Meru is today popularly called the Pascal's triangle \cite{coolidge1949, edwards2002, edwards2013, smith2010} after Blaise Pascal (1623-1662) who introduced this triangle in his treatise {\it Trait\'{e} du triangle arithm\'{e}tique} (1653)
\footnote{Keeping up with the tradition of giving due credit to the original propounder, the Pascal's triangle should be called the Pingala's triangle or the Pingala-Pascal's triangle.}. 
The notion of Pascal's triangle and its properties were also known to the Chinese--Jia
Xian (1010-1070) and Yang Hui (1238-1298), the Germans--Petrus Apianus (1495-1552) and Michael Stifel (1487-1567), and the Italian mathematicians Niccolo Fontana Tartaglia (1499-1557) and Gerolamo Cardano (1501-1576).

\subsection{Other additive sequences}
Lucas sequence, like Fibonacci sequence, is given by $l_n = l_{n-1} + l_{n-2}$ with $l_1 = 2$ and $l_2 = 1$ \cite{hoggatt1969, koshy2001}. In general, starting with \(g_1=a\) and \(g_2=b\), one can construct the following sequence: $a,~b,~a+b,~a+2b,~2a+3b,~3a+5b,~\cdots,~\mathbf{g_{n \ge 3}=f_{n-2}a+f_{n-1}b},~\cdots$. This general sequence is customarily called the Gopala-Hemchandra sequence \cite{singh1936, singh1985}. Furthermore, Narayana Pandita in his book {\it Ganita Kaumudi} (1356) \cite{ganitkaumudi} studies additive sequences where each term is the sum of the \(p\)-previous terms. He states the problem as: {\it A cow gives birth to a calf every year. The calves become young and they begin giving birth to calves when they are three years old. Tell me, O learned man, the number of progeny produced during twenty years by one cow.}

\subsection{Golden ratio}
The golden ratio \cite{huntley1970, runion1972, herz1987a, herz1987b, vajda1989, stakhov1989, livio2002, heath2002, kim2008}
\footnote{The golden ratio is also called {\it golden proportion}, {\it golden number}, {\it golden section}, {\it golden mean}, {\it divine proportion}, and {\it extreme and mean ratio}.}
as defined by Euclid in his book {\it The Elements} \cite{heath2002} is: 
{\it A straight line is said to have been cut in extreme and mean ratio
when, as the whole line is to the greater segment, so is the greater
to the less.} 

\begin{figure}%
\centering
\includegraphics[width = 3in]{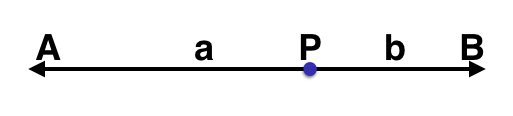}%
\caption{Division of a line into \(2\) segments.}
\label{fig-goldenratio2}
\end{figure} 

That is, the {\it golden ratio} arises when we consider division of a line segment $AB$ with a point $P$ such that  $\frac{AP}{BP} = \frac{AB}{AP}$, where $AP > BP$
(see Fig. \ref{fig-goldenratio2}). Given $AP = a$ and $BP = b$ are two positive numbers, the above problem translates as
\be \label{eq-2ratio} \frac{a}{b} = \frac{a+b}{a}. \ee 
Taking $\frac{a}{b} = x$, the above equation can be rewritten as $x = 1 + \frac{1}{x}$. This reduces to the characteristic equation 
\footnote{The {\it characteristic equation} is the {\it minimal polynomial} from which all the algebraic properties of an algebraic number (here $\Phi$) can be drawn.}
\be \label{eq-2ce} X(x) = x^2 - x -1=0, \ee 
whose positive solution is 
\be \Phi = \frac{\sqrt{5}+1}{2} = 1.618. \ee

\subsection{Relation between Fibonacci sequence and golden ratio} 
The Fibonacci sequence is closely related to the Golden Ratio in the sense that the limiting ratio value of the Fibonacci sequence, i.e., the ratio of successive numbers of the Fibonacci sequence tends to the golden ratio,
\be \lim_{n \rightarrow \infty} \frac{f_{n+1}}{f_n} = \Phi. \ee

\subsection{Properties of Fibonacci numbers and golden ratio}
The golden ratio has a number of interesting properties. They are listed below:

\begin{enumerate}
\item Some relations between Fibonacci numbers $f_0 = 0$, $f_1 = 1$, $f_{n \ge 2} = f_{n-1} + f_{n-2}$, and Lucas numbers $l_0 = 2$, $l_1 = 1$, $l_{n \ge 2} = l_{n-1} + l_{n-2}$.
\begin{enumerate}
\item $l_n = f_{n+1} + f_{n-1} = 2f_{n+1} - f_{n}$.

\item $f_n + f_{n+2} = l_{n+1}$.

\item $l_n + l_{n+2} = 5f_{n+1}$.
\end{enumerate}

\item $f_{n+1} f_{n-1} - f_n^2 = (-1)^n$ (Cassini's identity).

\item $\Phi^2 = \Phi + 1$.

\item $\Phi = \sqrt{1 + \sqrt{1 + \sqrt{1 + \cdots}}}$.

\item $\Phi = \frac{1 + \sqrt{5}}{2}$ and $\phi = - \frac{1 - \sqrt{5}}{2}$.

\item $\phi = \frac{1}{\Phi} = \Phi - 1$.

\item $\Phi = 1 + \frac{1}{\Phi}$.

\item The golden ratio, continued fractions and its convergents.
\begin{enumerate}
\item The continued fraction
\footnote{A {\it continued fraction} is a form of representing a number by nested fractions,
all of whose numerators are \(1\). 
The continued fraction of a rational number \(x\) is finite and is represented as
$x = a_0 + \frac{1}{a_1 + \frac{1}{a_2 + \frac{1}{\ddots + \frac{1}{a_n}}}} \equiv [a_0;a_1,a_2,\cdots, a_n]$, where $a_1,~a_2,~\cdots, ~a_n$ are positive integers and $a_0$ is any integer.
For example, $\frac{5}{3} = 1 + \frac{1}{1 + \frac{1}{2}} \equiv [1;1,2]$ and $\frac{10}{7} = 1 + \frac{1}{2 + \frac{1}{3}} \equiv [1;2,3]$. 
Note that the first term is followed by a semicolon, while other terms are followed by commas.
If \(x\) is irrational, then $n \rightarrow \infty$.}
of the golden ratio:
\be \Phi = 1 + \frac{1}{1 + \frac{1}{1 + \frac{1}{1 + \ddots}}} \equiv [1;\overline{1}]. \ee

\item The convergents
\footnote{A {\it convergent} is the truncation of a continued fraction. For example, the second convergent of $[1;2,3]$ is $[1;2]$ and the $m^{th}$ convergent of $[a_0;a_1,a_2,\cdots, a_n]$ is $[a_0;a_1,a_2,\cdots, a_{m-1}]$. That is, $[a_0;a_1,a_2,\cdots, a_n]_m := [a_0;a_1,a_2,\cdots, a_{m-1}]$.}
of the golden ratio:
\be [\Phi]_n = [1;\overline{1}]_n = \frac{f_{n+1}}{f_n}. \ee

\item The continued fraction of powers of the golden ratio:
\be  
[\Phi^n] = \left\{
  \begin{array}{ll}
   {[}l_n;\overline{l_n}{]} \quad &\text{(\(n\) odd)},\\
   {[}l_{n} - 1;\overline{1,l_{n} - 2}{]} \quad &\text{(\(n\) even)}.
  \end{array}
\right.
\ee

\item The convergents of powers of the golden ratio:
\be 
\frac{f_{a(n+1)}}{f_{an}} = \left\{
  \begin{array}{ll}
   {[}\Phi^a{]}_{n} \quad &\text{(\(a\) odd)},\\
   {[}\Phi^a{]}_{2n} \quad &\text{(\(a\) even)}.
  \end{array}
\right.
\ee
\end{enumerate}

\item $\Phi = \lim_{n \rightarrow \infty} \frac{f_{n+1}}{f_n}$.

\item $\Phi^{n} = \Phi^{n-1} + \Phi^{n-2} = \Phi  f_n + f_{n-1}$.

\item $f_n = \frac{\Phi^n - (-\phi)^{n}}{\sqrt{5}}$ (Binet's formula).

\item $\Phi^n = \frac{l_n + f_n \sqrt{5}}{2}$.

\item $\Phi$ as an infinite series: 
$\Phi = \frac{13}{8} + \sum_{n=0}^{\infty} \frac{(-1)^{n+1} (2n+1)!}{n!(n+2)!4^{2n+3}}$.

\item $\Phi$ as trigonometric functions:
$\Phi = 1 + 2\sin 18^{\circ} = 2\sin 54^{\circ} = \frac12\csc 18^{\circ}$.
\end{enumerate}

\subsection{Applications of golden ratio}
The golden ratio is certainly a famous number, and also a {\it divine} one as considered by some \cite{huntley1970, livio2002}. It allegedly appears everywhere. 

\begin{itemize}
\item {\it In geometry, maths and science} \cite{livio2002, kepler1966, doczi1986, osborne1986, srinivasan1992, stakhov1998, bradley1999, ball2003, kapusta2004, heyrovska2005, marek2006, sigalotti2006, sen2008, shao2008, hretcanu2009, esmaeili2010, daubechies2010, tanackov2011, agaian2011, agaian2012, mishra2012, schretter2012, hassaballah2013, yang2013, liang2016}. The golden ratio appears, by construction, in geometrical objects such as the golden polygons (triangle, rectangle, pentagon, etc.) and golden spirals. It also appears in science, physical theories and problems. 

\item {\it In nature} \cite{livio2002, brousseau1969, douady1996}. It exhibits in natural flora and fauna in the form of golden shapes such as spirals and pentagon. A tantalizing connection appears between the Fibonacci numbers (and hence the golden golden) and phyllotaxis (i.e., the arrangement of leaves on a stem, scales on a pine cone, florets on a sunflower, infloresences on a cauliflower, etc.). The plant tendrils get twisted by spirals, the helical motions are seen in the growth of roots and sprouts, and the sunflower seeds are arranged along the spirals. The spirals in sunflowers appear to rotate both clockwise (21 spirals) and counterclockwise (34 spirals). Remarkably, the numbers 21 and 34 are consecutive Fibonacci numbers. The golden ratio often shows in horns of rams, goats and antelopes. The pentagonal symmetry in the form of star fish, five-petal flowers, certain cactus plants, etc. are widespread in nature. 

\item {\it In human body} \cite{livio2002, chan2009, ajluni2010, kleider2010, henein2011}. It is believed that human body and its parts appear in the golden ratio. A ratio of feet-to-head height to feet-to-navel height (and also, ratio of feet-to-navel height to navel-to-head height) is called the navel ratio. A perfect human body is divided by the navel into the golden section. The human hand and face are also based on the golden ratio.

\item {\it In architecture} \cite{livio2002, herz2000, velmurugan2020}. It appears that the golden ratio has been used significantly in architecture: in Parthenon, in Great Pyramids of Egypt, in Indian meditation symbol {\it Sri Yantra}, in Taj Mahal and several ancient Indian temples such as Tanjavur Brihadeeshwara temple.    

\item {\it In art, painting and music} \cite{wilson1993, wilson2001, livio2002, frayling1992, zanten1999, olariu1999, jensen2002, kak2004, lian2013, kazlacheva2016}. The golden ratio is also prevalent in art, music and painting. For example, in the works of Da Vinci ({\it The Annunciation, Madonna with Child and Saints, The Mona Lisa, St. Jerome, An Old Man, and The Vitruvian Man}), in {\it The Holy Family} by Michelangelo, {\it The Crucifixion} by Raphael and {\it The Sacrament of the Last Supper} by Salvador Dali.  
It is illustrated in prolific number in portraits, paintings of Christian God and sculptures during the renaissance epoch. In music, it is present in works of Beethoven, Mozart, Wilson's {\it Meru 1} etc.
\end{itemize}
As asserted by many, it exists in any place where life and beauty are present.

\subsection{Are Fibonacci sequence and golden ratio sacred?}
Despite all-round great appearance of the golden ratio, many hold skeptical views on this \cite{devlin, falbo1, falbo2, gould1981, markowsky1992, spira2012}. The reasons are multifold.

\begin{itemize}
\item Application of the golden ratio to aesthetics is, by its nature, subjective and controversial. In order to find the golden ratio in our everyday life, we consider the following either separately or in combination \cite{spira2012}: (i) arbitrary placement of points, lines, rectangles and spirals,
(ii) arbitrary thickness of points and lines used as basis for measurements, and
(iii) measurements of monuments eroded by time and of objects in photographs distorted
by perspective.

\item Not all spirals in the nature are the golden ones. The nautilus shell, a prime pedagogical example, corresponds to a spiral with the value $\Phi' = 1.33~(<1.618)$.
\end{itemize}

\section{Early generalizations of golden ratio}
There have been sincere attempts to extend or generalize the notion of golden ratio from various perspectives such as generalizations of Euclid's problem, limits of recurrence relations, and the characteristic equations \cite{wilson1993, wilson2001, herz1987a, herz1987b, stakhov1989, fowler1982, kapur1988, engstrom1988, spinadel1999a, spinadel1999b, spinadel1999c, spinadel2000, spinadel2002, lipovetsky2000, rakocevic2004, krcadinac2006, stakhov2007, hashemiparast2011, dutta2020}. Fowler \cite{fowler1982} revisited the Euclid's problem {\it the line divided in extreme and mean ratio} and explored the propositions not investigated and proved in Euclid's {\it Elements}. Here, we review briefly some early generalizations.

\subsection{Golden p-proportions of Alexey Stakhov} 
Recall the Euclid's division problem of a line segment $AB$ into two segments $AP (= a)$ and $BP (= b)$ where $AP > BP$ (see Fig. \ref{fig-goldenratio2}). 
Alexey Stakhov, a Russian mathematician, considered the following generalization in his book \cite{stakhov1989}
\be \label{eq-stakhov1}
\frac{AP}{BP} = \big(\frac{AB}{AP} \big)^p \Rightarrow \frac{a}{b} = \big(\frac{a+b}{a} \big)^p,  \ee
where \(p\) is a non-negative integer. 
From Eq. (\ref{eq-stakhov1}), with $\frac{AB}{AP} = x$, we obtain the following algebraic equation
\be \label{eq-stakhov2} x^{p+1} = x^p + 1, \ee
whose the only positive solution $\chi_p$ is called the {\it golden p-proportion}. 
The {\it Fibonacci p-numbers} are obtained with the recurrence relation(s), 
\bea 
f_{n}(p) &=& t_{n-1}(p) + t_{n-(p+1)}(p),~~~(n \ge p+1) \nonumber \\ 
f_{n+1}(p) &=& t_{n}(p) + t_{n-p}(p),~~~(n \ge p) \label{eq-stakhov3}
\eea
where $f_k(p) = 1,~k =0,1,\cdots, p$. These numbers are related to the concept of {\it ``deformed'' Pascal's p-triangles} via the binomial coefficients as
\be \label{eq-stakhov4} f_{n+1}(p) = \sum_{k=0}^{\infty} \binom{n-kp}{k}. \ee

\begin{table}[]
\centering
\begin{tabular}{|l|l|l|l|l|l|l|l|l|l|l|}
\hline
$n$        & 0 & 1 & 2 & 3 & 4  & 5  & 6  & 7   & 8   & 9   \\ \hline
$f_n(p=0)$ & 1 & 2 & 4 & 8 & 16 & 32 & 64 & 128 & 256 & 512 \\ \hline
$f_n(p=1)$ & 1 & 1 & 2 & 3 & 5  & 8  & 13 & 21  & 34  & 55  \\ \hline
$f_n(p=2)$ & 1 & 1 & 1 & 2 & 3  & 4  & 6  & 9   & 13  & 19  \\ \hline
$f_n(p=3)$ & 1 & 1 & 1 & 1 & 2  & 3  & 4  & 5   & 7   & 10  \\ \hline
\end{tabular}%
\caption{The Fibonacci \(p\)-numbers $f_n(p)$ for different \(p\) values. $f_n(p=0) = 2^n$ are the binary numbers, $f_n(p=1)$ are the Fibonacci numbers, and so on.}
\label{table-stakhov}
\end{table}

Note the following observations:
\begin{enumerate}
\item $f_n(p=0) = 2^n$ are the binary numbers, $f_n(p=1)$ are the Fibonacci numbers, and so on (see Table \ref{table-stakhov}).

\item $\chi_{0} = 2$, $\chi_1 = \frac{1 + \sqrt{5}}{2} = \Phi$, $\chi_{\infty} =1$, and $1 \le \chi_p \le 2$. 

\item $\chi_{p}^n = \chi_{p}^{n-1} + \chi_{p}^{n-(p+1)} = \chi_{p} \times \chi_{p}^{n-1}$.

\item Binomial coefficients, Fibonacci p-numbers, and golden p-proportions.\\
$f_{n+1}(p) = \sum_{k=0}^{\infty} \binom{n-kp}{k}$ and $\chi_p = \lim_{n \rightarrow \infty} \frac{f_{n+1}(p)}{f_{n}(p)}$.
\end{enumerate}

The notions of the golden \(p\)-proportions and Fibonacci \(p\)-numbers generalized the original mathematical concepts, and led to several interesting applications including in the different fields of mathematics and computer science \cite{stakhov1989, stakhov1998, stakhov2007}.

\subsection{Metallic means family of Vera Spinadel} 
Vera Spinadel, an Argentinean mathematician, considered an interesting generalization of the Fibonacci recurrence relation, $t_{n+1} = t_{n} + t_{n-1}$, in the following form
\bea \label{eq-spinadel1}
t_{n+1} &=& pt_{n} + qt_{n-1}, \\
\Rightarrow \frac{t_{n+1}}{t_{n}} &=& p + q \frac{t_{n-1}}{t_{n}}, \nonumber
\eea
where \(p\) and \(q\) are non-negative integers. Assuming that $\lim_{n \rightarrow \infty} \frac{t_{n+1}}{t_{n}} = x$ exists, we have
\be \label{eq-spinadel2}
x = p + \frac{q}{x} \Rightarrow x^2 = px + q. \ee
The algebraic equation (\ref{eq-spinadel2}) has a solution
\be \label{eq-spinadel3} \chi_{p,q} = \frac{p + \sqrt{p^2 + 4q}}{2}. \ee 
Positive solutions in Eq. (\ref{eq-spinadel3}) form a {\it metallic means family} (MMF), and Vera Spinadel gave a number of applications of the metallic means in her works \cite{spinadel1999a, spinadel1999b, spinadel1999c, spinadel2000, spinadel2002}. 
Note the following observations:
\begin{enumerate}
\item $x^2 = px +q$ implies $x = \sqrt{q+ p\sqrt{q+ p\sqrt{q+ p\sqrt{q+ \cdots}}}}$.

\item $x = p + \frac{q}{x}$ implies $x = p + \frac{q}{p + \frac{q}{p + \frac{q}{p + \ddots}}}$.

\item $\chi_{p,q} = \frac{p + \sqrt{p^2 + 4q}}{2} = \lim_{n \rightarrow \infty} \frac{t_{n+1}}{t_{n}}$.

\item $\chi_{1,1} = \frac{1 + \sqrt{5}}{2} = [1;\overline{1}] = \Phi$.

\item $\chi_{p,1} = \frac{p + \sqrt{p^2 + 4}}{2} = [p;\overline{p}]$.

\item $\chi_{4,1} = 2 + \sqrt{5} = [4;\overline{4}] = \Phi^3$.
\end{enumerate}

\subsection{Mount merus of Erwin Wilson} 
Recall Pingala's Mount Meru (Pascal's Triangle) in Fig. \ref{fig-pingala}. It was illustrated in 1968 by Thomas Green \cite{green1968} that the sum of the simplest diagonals of Mount Meru yields the Fibonacci sequence, and that the sum of other diagonals similarly generate other recurrence relations, each with its own limit. Ervin Wilson, a Mexican/American music theorist, investigated these other diagonals and their recurrence relations in music \cite{wilson1993, wilson2001}. He considered recurrence relations that he called {\it Meru 1} through {\it Meru 9}. See Table \ref{table-wilson}.

\begin{table}[]
\centering
\begin{tabular}{|l|l|l|l|}
\hline
meru   & recurrence relation       & characteristic equation & convergence limit \\ \hline
Meru 1 & $A_n = A_{n-1} + A_{n-2}$ & $x^2 = x + 1$           & 1.61803           \\ \hline
Meru 2 & $B_n = B_{n-1} + B_{n-3}$ & $x^3 = x^2 + 1$         & 1.46557           \\ \hline
Meru 3 & $C_n = C_{n-2} + C_{n-3}$ & $x^3 = x + 1$           & 1.32472           \\ \hline
Meru 4 & $D_n = D_{n-1} + D_{n-4}$ & $x^4 = x^3 + 1$         & 1.38028           \\ \hline
Meru 5 & $E_n = E_{n-3} + E_{n-4}$ & $x^4 = x + 1$           & 1.22074           \\ \hline
Meru 6 & $F_n = F_{n-1} + F_{n-5}$ & $x^5 = x^4 + 1$         & 1.32472           \\ \hline
Meru 7 & $G_n = G_{n-2} + G_{n-5}$ & $x^5 = x^3 + 1$         & 1.23651           \\ \hline
Meru 8 & $H_n = H_{n-3} + H_{n-5}$ & $x^5 = x^2 + 1$         & 1.19386           \\ \hline
Meru 9 & $I_n = I_{n-4} + I_{n-5}$ & $x^5 = x + 1$           & 1.16730           \\ \hline
\end{tabular}%
\caption{Wilson's {\it Meru 1} through {\it Meru 9}, each with its recurrence relation, the characteristic equation, and the convergence limit.}
\label{table-wilson}
\end{table}

\subsection{Lower and upper golden ratios of Vedran Krcadinac}
We saw earlier that Stakhov considered generalization of the form $\frac{a}{b} = \big(\frac{a+b}{a} \big)^p$ leading to the algebraic equation $x^{p+1} - x^p - 1 = 0$ when $\frac{a+b}{a} = x$. A similar generalization, proposed by Krcadinac \cite{krcadinac2006}, for non-negative integer \(p\), is 
\be \label{eq-krcadinac1}
\big(\frac{a}{b} \big)^p= \frac{a+b}{a}. \ee
This relation, in general, leads to two algebraic equations:
\bea 
X_1(x) &=& x^{p+1} - x -1 = 0,~~(\text{when}~\frac{a}{b} = x) \label{eq-krcadinac2} \\ 
X_2(x) &=& x(x - 1)^{p} -1 = 0,~~(\text{when}~\frac{a+b}{a} = x) \label{eq-krcadinac3}
\eea
Let $\varphi_p$ be the positive root of the polynomial $X_1(x)$ and $\phi_p$ be that of the polynomial $X_2(x)$. Then, $\varphi_p$ and $\phi_p$ are respectively called the $p^{th}$ {\it lower} and {\it upper} golden ratios. Note the following observations:
\begin{enumerate}
\item $\varphi_0 = \textit{undefined}$ and $\phi_0 = 1$. 

\item $\lim_{p \rightarrow \infty} \varphi_p = 1$ and $\lim_{p \rightarrow \infty} \phi_p = 2$. 

\item Evidently, $(\varphi_p)^p = \phi_p$.

\item Recurrence relation for $X_1(x)$: $f_{n}(p) = f_{n-p}(p) + f_{n-(p+1)}(p)$.

\item Recurrence relation for $X_2(x)$: $F_{n}(p) = \sum_{k=1}^p \binom{p}{k}~ (-1)^{k+1} F_{n-k}(p) + F_{n-(p+1)}(p)$.

\item $\lim_{n \rightarrow \infty} \frac{f_{n+1}(p)}{f_{n}(p)}= \varphi_p$ and $\lim_{n \rightarrow \infty} \frac{F_{n+1}(p)}{F_{n}(p)}= \phi_p$.

\item $\lim_{n \rightarrow \infty} \frac{f_{n+p}(p)}{f_{n}(p)}= (\varphi_p)^p = \phi_p$.
\end{enumerate}

\section{\(p\)-sequences}
We call the Fibonacci sequence a \(2\)-sequence because it is generated by the sum of two previous terms. In a similar spirit, we introduce the \(p\)-sequence
\footnote{\(p\) in the \(p\)-sequence is for {\it Pingala}, {\it Phi($\Phi$)}, and {\it previous}.}.\\ 

To construct a \(p\)-sequence, we begin with \(p\) seeds $(s_0,s_1,\cdots,s_{p-1})$ such that $t_0=s_0,~t_1=s_1,\cdots,t_{p-1}=s_{p-1}$, and the $n^{th}$ term is the sum of its {\it \(p\) previous} terms
\footnote{This can be equivalently rewritten as
$t_{n+p}(p) := t_{n+p-1}(p) + t_{n+p-2}(p) + \cdots + t_{n}(p).$ 
}:
\be t_n(p) := t_{n-1}(p) + t_{n-2}(p) + \cdots + t_{n-p}(p) = \sum_{k=n-p}^{n-1} t_k(p). \ee 
By definition of $t_n(p)$, we have
\bea
t_{n+1}(p) &>&  t_n(p), \\
t_{n+1}(p) &=& 2t_n(p) -t_{n-p}(p) <  2t_n(p).
\eea
Depending on the values of seeds, one can construct an infinite number of \(p\)-sequences. A few typical \(p\)-sequences are:\\

(i) {\it General} \(p\)-sequence whose seeds are arbitrary.
\be S_G(p) \equiv \{(s_0,s_1,\cdots, s_{p-1}),~t_n(p)\}. \ee

(ii) $\textit{k}$ \(p\)-sequence whose $k^{th}$ seed is unity and other seeds are zero. 
\be S_k(p) \equiv \{(s_i = \delta_{ik},~ 0 \le i \le p-1),~t_n(p)\}. \ee
For example, $S_0(p) \equiv \{(1,0,\cdots,0),~t_n(p)\}$, $S_1(p) \equiv \{(0,1,0,\cdots,0),~t_n(p)\}$, and $S_{p-1}(p) \equiv \{(0,0,\cdots,1),~t_n(p)\}$.
Interestingly, we can rewrite $t_n[S_G(p)]$ in terms of seeds using these $\textit{k}$ \(p\)-sequences,
\be t_n[S_G(p)] = \sum_{k=0}^{p-1} t_n[S_k(p)]s_k~~(n \ge 0). \ee   
For example, $t_1[S_G(p)] = 0.s_0 + 1.s_1 + \cdots + 0.s_{p-1} = s_1$.\\

(iii)  {\it Coefficient} \(p\)-sequence whose all seeds are unity.
\be S_C(p) \equiv \{(s_k = 1,~ 0\le k \le p-1),~t_n(p)\}. \ee
There is an important relation between the terms of {\it coefficient} \(p\)-sequence and those of $\textit{k}$ \(p\)-sequences: $S_C(p) \equiv \sum_{k=1}^p S_k(p)$. Put differently,
\be t_n[S_C(p)] = \sum_{k=0}^{p-1} t_n[S_k(p)]. \ee

(iv) {\it Exponent} \(p\)-sequence whose seeds are $(0,1,\cdots,p-1)$.
\be S_X(p) \equiv \{(s_k = k,~ 0\le k \le p-1),~t_n(p)\}. \ee

(v) {\it Syllable} \(p\)-sequence whose seeds are $(1,2,\cdots,2^{p-1})$.
\be S_S(p) \equiv \{(s_k = 2^{k},~ 0\le k \le p-1),~t_n(p)\}. \ee
We will learn the significance of these particular sequences in the forthcoming articles.
For illustrations of and getting familiarized with these sequences, see Tables \ref{table-p2}, \ref{table-p3}, \ref{table-p4} and \ref{table-p5}. 
Henceforth, $S_G(p) \equiv S(Gp)$, $S_k(p) \equiv S(kp)$, $S_C(p) \equiv S(Cp)$, $S_X(p) \equiv S(Xp)$, and $S_S(p) \equiv S(Sp)$ will be used interchangeably.
We will denote the $n^{th}$-term of an arbitrary \(p\)-sequence by $t_n(p)$, and that of a particular \(p\)-sequence, viz. {\it exponent} sequence by $t_n(Xp)$, {\it syllable} sequence by $t_n(Sp)$, 
{\it k} sequence by $t_n(kp)$, and so on.
Furthermore, in case of no ambiguity, we will not mention \(p\) explicitly in the sequence names and their terms. 

\subsection{\(1\)-sequence}
What is \(1\)-sequence? 
We construct a \(1\)-sequence by choosing a seed $s_0 \ge 0$ and a constant $a \ge 0$ such that $t_0 = s_0$, and for $n \ge 1$
\bea 
t_1 &=& t_0 + a = s_0 + a, \nonumber \\
t_2 &=& t_1 + a = s_0 + 2a, \nonumber \\
t_n &=& t_{n-1} +a = s_0 + na. 
\eea 
Thus, an additive \(1\)-sequence is essentially an arithmetic progression. With $s_0 = 0$ and $a =1$, \(1\)-sequence is the set of whole numbers 
\be S(1) = \{0,1,2,\cdots, 99,100, \cdots \}. \ee
When $a =0$, the \(1\)-sequence is a {\it constant} sequence: $\{s_0, s_0, s_0, \cdots\}$.

\section{Generating functions of \(p\)-sequences}
The generating function for \(p\)-sequences can be given by the power series
\be f_{p}(x) = \sum_{n=0}^{\infty} t_n(p) x^n, \ee
where $t_n(p)$ is the $n^{th}$ term of a given \(p\)-sequence.
If we assume that the power series converges, we can show that $f_{p}(x)$ is given by
\be \big(1 - \sum_{k=1}^p x^k \big)f_{p}(x) = \sum_{k=0}^{p-1} \big[t_k(p) - \sum_{j=0}^{k-1} t_j(p) \big]x^k. \ee

For example, for the {\it exponent} \(p\)-sequence $S(Xp)$, the generating functions are
\bea
f_{X2}(x) &=& \frac{x}{1-x-x^2}, \nonumber \\
f_{X3}(x) &=& \frac{x+x^2}{1-x-x^2-x^3}, \nonumber \\
f_{X4}(x) &=& \frac{x+x^2}{1-x-x^2-x^3-x^4}, \nonumber \\
f_{X5}(x) &=& \frac{x+x^2-2x^4}{1-x-x^2-x^3-x^4-x^5}. \nonumber
\eea

\section{Limiting ratio value of \(p\)-sequences}
For an arbitrary \(p\)-sequence whose subsequent terms are the sum of \(p\)-previous terms [$t_n(p) := \sum_{k=n-p}^{n-1} t_k$], we see from Tables \ref{table-p2}, \ref{table-p3}, \ref{table-p4} and \ref{table-p5} that the limiting ratio value of every \(p\)-sequence approaches a constant, say $\Phi_p$. That is,
\be \label{eq-p-ratio-limit}
\Phi_p = \lim_{n \rightarrow \infty} \frac{t_{n+1}(p)}{t_{n}(p)}. \ee

Because $t_{n+1}(p) >  t_n(p)$ and $t_{n+1}(p) < 2t_n(p)$, hence
\be 1 < \Phi_p < 2. \ee

Using Eq. (\ref{eq-p-ratio-limit}), for integers \(u\) and \(v\), it is easy to see that 
\bea \label{eq-p-ratio-limit1}
\lim_{n \rightarrow \infty} \frac{t_{n+u}(p)}{t_{n}(p)} &=& \Phi_p^{u},~~~ 
\lim_{n \rightarrow \infty} \frac{t_{n+u}(p)}{t_{n+v}(p)} = \Phi_p^{u-v}. 
\eea

\subsection{$\Phi_p$ in the limit $p \rightarrow \infty$}
Consider the {\it syllable} \(p\)-sequence $S_S(p \rightarrow \infty) = \{1,~2,~4,~8,~16, \cdots \}$, 
and the {\it extended syllable} \(p\)-sequence
$S(p \rightarrow \infty) \equiv \{(0,~1),~S_S(p \rightarrow \infty) \} = \{0,~1,~1,~2,~4,~8,~16, \cdots \},$
where each term is the sum of all the previous terms except the first two. For both these sequences, we have
\be \Phi_{p \rightarrow \infty} = \lim_{n \rightarrow \infty} \frac{t_{n+1}(p)}{t_{n}(p)} = 2. \ee 
Moreover, from Table \ref{table-plimits}, $\Phi_p =2$ for \(p \ge 18\).

\subsection{Propositions}
In Tables \ref{table-p2}, \ref{table-p3}, \ref{table-p4} and \ref{table-p5}, we constructed \(p\)-sequences for $p=2,~3,~4,~5$, and found their limiting ratio values.  
Similarly, one can construct tables of higher \(p\)-sequences and find their limiting ratio values.
See Table \ref{table-plimits} for the values of limiting ratios. In this regard, we propound the following two propositions.\\

(P1) The limiting ratio value of any \(p\)-sequence $S(p)$ is $\Phi_p$. It is independent of the initial conditions (i.e., the seeds).\\

(P2) $\Phi_{p \ge 18} =2.$

\begin{table}[]
\centering
\begin{tabular}{|l|l|l|l|l|l|}
\hline
$n$     & $S_1(2)$     & $S_0(2)$    & $S_C(2)$    & $S_S(2)$    & $S_G(2)$   \\ \hline 
\hline
0  & 0     & 1     & 1      & 1      & 2       \\ \hline
1  & 1     & 0     & 1      & 2      & 21      \\ \hline 
\hline
2  & 1     & 1     & 2      & 3      & 23      \\ \hline 
\hline
3  & 2     & 1     & 3      & 5      & 44      \\ \hline
4  & 3     & 2     & 5      & 8      & 67      \\ \hline
5  & 5     & 3     & 8      & 13     & 111     \\ \hline
6  & 8     & 5     & 13     & 21     & 178     \\ \hline
7  & 13    & 8     & 21     & 34     & 289     \\ \hline
8  & 21    & 13    & 34     & 55     & 467     \\ \hline
9  & 34    & 21    & 55     & 89     & 756     \\ \hline
10 & 55    & 34    & 89     & 144    & 1223    \\ \hline
11 & 89    & 55    & 144    & 233    & 1979    \\ \hline
12 & 144   & 89    & 233    & 377    & 3202    \\ \hline
13 & 233   & 144   & 377    & 610    & 5181    \\ \hline
14 & 377   & 233   & 610    & 987    & 8383    \\ \hline
15 & 610   & 377   & 987    & 1597   & 13564   \\ \hline
16 & 987   & 610   & 1597   & 2584   & 21947   \\ \hline
17 & 1597  & 987   & 2584   & 4181   & 35511   \\ \hline
18 & 2584  & 1597  & 4181   & 6765   & 57458   \\ \hline
19 & 4181  & 2584  & 6765   & 10946  & 92969   \\ \hline
20 & 6765  & 4181  & 10946  & 17711  & 150427  \\ \hline
21 & 10946 & 6765  & 17711  & 28657  & 243396  \\ \hline
22 & 17711 & 10946 & 28657  & 46368  & 393823  \\ \hline
23 & 28657 & 17711 & 46368  & 75025  & 637219  \\ \hline
24 & 46368 & 28657 & 75025  & 121393 & 1031042 \\ \hline
25 & 75025 & 46368 & 121393 & 196418 & 1668261 \\ \hline
\end{tabular}%
\caption{\(2\)-sequences. (i) $S_C \equiv S_1 + S_0$. (ii) $S_X = S_1$. (iii) $S_1 \sim S_0 \sim S_C \sim S_S$. (iv) $S_G$ is a general \(2\)-sequence with seeds $s_1=2,~s_2=21$. (v) For each of these \(2\)-sequences, $\lim_{n \rightarrow \infty} \frac{t_{n+1}}{t_{n}}=1.61803$.}
\label{table-p2}
\end{table}

\begin{table}[]
\centering
\begin{tabular}{|l|l|l|l|l|l|l|}
\hline
$n$     & $S_2(3)$     & $S_1(3)$     & $S_0(3)$    & $S_C(3)$      & $S_X(3)$      & $S_S(3)$       \\ \hline
\hline
0  & 0      & 0      & 0      & 1       & 0       & 1       \\ \hline
1  & 0      & 1      & 0      & 1       & 1       & 2       \\ \hline
2  & 1      & 0      & 1      & 1       & 2       & 4       \\ \hline 
\hline
3  & 1      & 1      & 1      & 3       & 3       & 7       \\ \hline 
\hline
4  & 2      & 2      & 2      & 5       & 6       & 13      \\ \hline
5  & 4      & 3      & 4      & 9       & 11      & 24      \\ \hline
6  & 7      & 6      & 7      & 17      & 20      & 44      \\ \hline
7  & 13     & 11     & 13     & 31      & 37      & 81      \\ \hline
8  & 24     & 20     & 24     & 57      & 68      & 149     \\ \hline
9  & 44     & 37     & 44     & 105     & 125     & 274     \\ \hline
10 & 81     & 68     & 81     & 193     & 230     & 504     \\ \hline
11 & 149    & 125    & 149    & 355     & 423     & 927     \\ \hline
12 & 274    & 230    & 274    & 653     & 778     & 1705    \\ \hline
13 & 504    & 423    & 504    & 1201    & 1431    & 3136    \\ \hline
14 & 927    & 778    & 927    & 2209    & 2632    & 5768    \\ \hline
15 & 1705   & 1431   & 1705   & 4063    & 4841    & 10609   \\ \hline
16 & 3136   & 2632   & 3136   & 7473    & 8904    & 19513   \\ \hline
17 & 5768   & 4841   & 5768   & 13745   & 16377   & 35890   \\ \hline
18 & 10609  & 8904   & 10609  & 25281   & 30122   & 66012   \\ \hline
19 & 19513  & 16377  & 19513  & 46499   & 55403   & 121415  \\ \hline
20 & 35890  & 30122  & 35890  & 85525   & 101902  & 223317  \\ \hline
21 & 66012  & 55403  & 66012  & 157305  & 187427  & 410744  \\ \hline
22 & 121415 & 101902 & 121415 & 289329  & 344732  & 755476  \\ \hline
23 & 223317 & 187427 & 223317 & 532159  & 634061  & 1389537 \\ \hline
24 & 410744 & 344732 & 410744 & 978793  & 1166220 & 2555757 \\ \hline
25 & 755476 & 634061 & 755476 & 1800281 & 2145013 & 4700770 \\ \hline
\end{tabular}%
\caption{\(3\)-sequences. (i) $S_C \equiv S_2 + S_1 + S_0$. (ii) $S_2 \sim S_0  \sim S_S$. (iii) $S_1 \sim S_X$. (iv) For each of these \(3\)-sequences, $\lim_{n \rightarrow \infty} \frac{t_{n+1}}{t_{n}}=1.83929$.}
\label{table-p3}
\end{table}

\begin{table}[]
\centering
\begin{tabular}{|l|l|l|l|l|l|l|l|}
\hline
$n$     & $S_3(4)$     & $S_2(4)$     & $S_1(4)$     & $S_0(4)$    & $S_C(4)$      & $S_X(4)$      & $S_S(4)$       \\ \hline
\hline
0  & 0       & 0      & 0      & 1      & 1       & 0       & 1        \\ \hline
1  & 0       & 0      & 1      & 0      & 1       & 1       & 2        \\ \hline
2  & 0       & 1      & 0      & 0      & 1       & 2       & 4        \\ \hline
3  & 1       & 0      & 0      & 0      & 1       & 3       & 8        \\ \hline 
\hline
4  & 1       & 1      & 1      & 1      & 4       & 6       & 15       \\ \hline 
\hline
5  & 2       & 2      & 2      & 1      & 7       & 12      & 29       \\ \hline
6  & 4       & 4      & 3      & 2      & 13      & 23      & 56       \\ \hline
7  & 8       & 7      & 6      & 4      & 25      & 44      & 108      \\ \hline
8  & 15      & 14     & 12     & 8      & 49      & 85      & 208      \\ \hline
9  & 29      & 27     & 23     & 15     & 94      & 164     & 401      \\ \hline
10 & 56      & 52     & 44     & 29     & 181     & 316     & 773      \\ \hline
11 & 108     & 100    & 85     & 56     & 349     & 609     & 1490     \\ \hline
12 & 208     & 193    & 164    & 108    & 673     & 1174    & 2872     \\ \hline
13 & 401     & 372    & 316    & 208    & 1297    & 2263    & 5536     \\ \hline
14 & 773     & 717    & 609    & 401    & 2500    & 4362    & 10671    \\ \hline
15 & 1490    & 1382   & 1174   & 773    & 4819    & 8408    & 20569    \\ \hline
16 & 2872    & 2664   & 2263   & 1490   & 9289    & 16207   & 39648    \\ \hline
17 & 5536    & 5135   & 4362   & 2872   & 17905   & 31240   & 76424    \\ \hline
18 & 10671   & 9898   & 8408   & 5536   & 34513   & 60217   & 147312   \\ \hline
19 & 20569   & 19079  & 16207  & 10671  & 66526   & 116072  & 283953   \\ \hline
20 & 39648   & 36776  & 31240  & 20569  & 128233  & 223736  & 547337   \\ \hline
21 & 76424   & 70888  & 60217  & 39648  & 247177  & 431265  & 1055026  \\ \hline
22 & 147312  & 136641 & 116072 & 76424  & 476449  & 831290  & 2033628  \\ \hline
23 & 283953  & 263384 & 223736 & 147312 & 918385  & 1592363 & 3919944  \\ \hline
24 & 547337  & 507689 & 431265 & 283953 & 1770244 & 3068654 & 7555935  \\ \hline
25 & 1055026 & 978602 & 831290 & 547337 & 3412255 & 5623572 & 14564533 \\ \hline
\end{tabular}%
\caption{\(4\)-sequences. (i) $S_C \equiv S_3 + S_2 + S_1 + S_0$. (ii) $S_3 \sim S_0  \sim S_S$. (iii) $S_1 \sim S_X$. (iv) For each of these \(4\)-sequences, $\lim_{n \rightarrow \infty} \frac{t_{n+1}}{t_{n}}=1.92756$.}
\label{table-p4}
\end{table}

\begin{table}[]
\centering
\begin{tabular}{|l|l|l|l|l|l|l|l|l|}
\hline
$n$  & $S_4(5)$     & $S_3(5)$     & $S_2(5)$     & $S_1(5)$     & $S_0(5)$    & $S_C(5)$      & $S_X(5)$      & $S_S(5)$       \\ \hline
\hline
0  & 0      & 0      & 0      & 0      & 1      & 1       & 0       & 1        \\ \hline
1  & 0      & 0      & 0      & 1      & 0      & 1       & 1       & 2        \\ \hline
2  & 0      & 0      & 1      & 0      & 0      & 1       & 2       & 4        \\ \hline
3  & 0      & 1      & 0      & 0      & 0      & 1       & 3       & 8        \\ \hline
4  & 1      & 0      & 0      & 0      & 0      & 1       & 4       & 16       \\ \hline
\hline
5  & 1      & 1      & 1      & 1      & 1      & 5       & 10      & 31       \\ \hline
\hline
6  & 2      & 2      & 2      & 2      & 1      & 9       & 20      & 61       \\ \hline
7  & 4      & 4      & 4      & 3      & 2      & 17      & 39      & 120      \\ \hline
8  & 8      & 8      & 7      & 6      & 4      & 33      & 76      & 236      \\ \hline
9  & 16     & 15     & 14     & 12     & 8      & 65      & 149     & 464      \\ \hline
10 & 31     & 30     & 28     & 24     & 16     & 129     & 294     & 912      \\ \hline
11 & 61     & 59     & 55     & 47     & 31     & 253     & 578     & 1793     \\ \hline
12 & 120    & 116    & 108    & 92     & 61     & 497     & 1136    & 3525     \\ \hline
13 & 236    & 228    & 212    & 181    & 120    & 977     & 2233    & 6930     \\ \hline
14 & 464    & 448    & 417    & 356    & 236    & 1921    & 4390    & 13624    \\ \hline
15 & 912    & 881    & 820    & 700    & 464    & 3777    & 8631    & 26784    \\ \hline
16 & 1793   & 1732   & 1612   & 1376   & 912    & 7425    & 16968   & 52656    \\ \hline
17 & 3525   & 3405   & 3169   & 2705   & 1793   & 14597   & 33358   & 103519   \\ \hline
18 & 6930   & 6694   & 6230   & 5318   & 3525   & 28697   & 65580   & 203513   \\ \hline
19 & 13624  & 13160  & 12248  & 10455  & 6930   & 56417   & 128927  & 400096   \\ \hline
20 & 26784  & 25872  & 24079  & 20554  & 13624  & 110913  & 253464  & 786568   \\ \hline
21 & 52656  & 50863  & 47338  & 40408  & 26784  & 218049  & 498297  & 1546352  \\ \hline
22 & 103519 & 99994  & 93064  & 79440  & 52656  & 428673  & 979626  & 3040048  \\ \hline
23 & 203513 & 196583 & 182959 & 156175 & 103519 & 842749  & 1925894 & 5976577  \\ \hline
24 & 400096 & 386472 & 359688 & 307032 & 203513 & 1656801 & 3786208 & 11749641 \\ \hline
25 & 786568 & 754784 & 707128 & 603609 & 400096 & 3257185 & 7443489 & 23099186 \\ \hline
\end{tabular}%
\caption{\(5\)-sequences. (i) $S_C \equiv S_4 + S_3 + S_2 + S_1 + S_0$. (ii) $S_4 \sim S_0  \sim S_S$. (iii) For each of these \(5\)-sequences, $\lim_{n \rightarrow \infty} \frac{t_{n+1}}{t_{n}}=1.96595$.}
\label{table-p5}
\end{table}

\begin{table}[]
\centering
\begin{tabular}{|l|l|l|l|l|}
\hline
$p$       & 2       & 3       & 4       & 5       \\ \hline
$\Phi_p$ & 1.61803 & 1.83929 & 1.92756 & 1.96595 \\ \hline
\hline
$p$       & 6       & 7       & 8       & 9       \\ \hline
$\Phi_p$ & 1.98358 & 1.99196 & 1.99603 & 1.99803 \\ \hline
\hline
$p$       & 10      & 11      & 12      & 13      \\ \hline
$\Phi_p$ & 1.99902 & 1.99951 & 1.99976 & 1.99988 \\ \hline
\hline
$p$       & 14      & 15      & 16      & 17      \\ \hline
$\Phi_p$ & 1.99994 & 1.99997 & 1.99998 & 1.99999 \\ \hline
\hline
$p$       & 18      & 19      & 20      & 21      \\ \hline
$\Phi_p$ & 2.0     & 2.0     & 2.0     & 2.0     \\ \hline
\end{tabular}%
\caption{The limitining ratio value $\Phi_p := \lim_{n \rightarrow \infty} \frac{t_{n+1}}{t_{n}}$ for \(p\)-sequences, $2 \le p \le 21$. We have limited ourselves here to five decimal places (for no sacred reasons). Evidently, $\Phi_2 < \Phi_3 < \cdots < \Phi_{p \ge 18} = 2$.}
\label{table-plimits}
\end{table}

\section{More additive sequences}
In this section, motivated by Pingala's syllable problem, Fibonacci's rabbit pair problem, and Narayan Pandit's cow's progeny problem, we investigate a general problem:
{\it A creature gives birth to $\alpha$ female young ones in one unit of time. Baby creature grows and gives birth when $\beta$ units of time old. The creature ceases to give birth after $\gamma$ terms, and dies when $\delta$ units of time old.. What is the total number of progeny at the end of \(n\) units of time? Initially, there is a single adult creature.} Following, we consider a few illustrations sans $\delta$. See Tables \ref{table-a1b2gna}, \ref{table-a1b3gna}, \ref{table-a2b2g3} and \ref{table-a2b3g2}. We invite the readers to investigate the problem taking into account $\delta$ also.

\begin{table}[]
\centering
\begin{tabular}{|l|l|l|l|}
\hline
$n$ & creature & baby $(b_n)$ at start & total $(t_n)$ at end\\ \hline
0   & 1        & 0    & 0     \\ \hline
\hline
1   & 1        & 1    & 1     \\ \hline
2   & 1        & 1    & 2     \\ \hline
3   & 1+1      & 2    & 3     \\ \hline
4   & 2+1      & 3    & 5     \\ \hline
5   & 3+2      & 5    & 8     \\ \hline
6   & 5+3      & 8    & 13    \\ \hline
7   & 8+5      & 13   & 21    \\ \hline
8   & 13+8     & 21   & 34    \\ \hline
9   & 21+13    & 34   & 55    \\ \hline
10  & 34+21    & 55   & 89    \\ \hline
11  & 55+34    & 89   & 144   \\ \hline
12  & 89+55    & 144  & 233   \\ \hline
13  & 144+89   & 233  & 377   \\ \hline
14  & 233+144  & 377  & 610   \\ \hline
15  & 377+233  & 610  & 987   \\ \hline
\end{tabular}%
\caption{$\alpha =1$, $\beta =2$ and $\gamma =\text{NA}$. This yields Pingala (Fibonacci) sequence. Here $t_{n \ge 1} = b_n + b_{n-1}$, $t_{n \ge 3} = b_n + b_{n-2} + b_{n-3}$ (sum of $\beta + 1$ terms), $t_{n \ge 3} = t_{n-1} + t_{n-2}$ (recurrence relation), $x^2 = x + 1$ (characteristic equation) and $x = \lim_{n \rightarrow \infty} \frac{t_{n+1}}{t_n} = 1.618$ (limiting ratio).}
\label{table-a1b2gna}
\end{table}

\begin{table}[]
\centering
\begin{tabular}{|l|l|l|l|}
\hline
$n$ & creature & baby $(b_n)$ & total $(t_n)$ \\ \hline
0   & 1        & 0    & 0     \\ \hline
\hline
1   & 1        & 1    & 1     \\ \hline
2   & 1        & 1    & 2     \\ \hline
3   & 1        & 1    & 3     \\ \hline
4   & 1        & 1    & 4     \\ \hline
5   & 1+1      & 2    & 5     \\ \hline
6   & 2+1      & 3    & 7     \\ \hline
7   & 3+1      & 4    & 10    \\ \hline
8   & 4+1      & 5    & 14    \\ \hline
9   & 5+2      & 7    & 19    \\ \hline
10  & 7+3      & 10   & 26    \\ \hline
11  & 10+4     & 14   & 36    \\ \hline
12  & 14+5     & 19   & 50    \\ \hline
\end{tabular}%
\caption{$\alpha =1$, $\beta =3$ and $\gamma =\text{NA}$. This sequence corresponds to Narayana Pandit's cow's progeny problem posed in {\it Ganit Kaumudi}. Here $t_{n \ge 3} = b_n + b_{n-1} + b_{n-2} + b_{n-3}$ (sum of $\beta + 1$ terms), $t_{n \ge 5} = t_{n-1} + t_{n-4}$ (recurrence relation), $x^4 = x^3 + 1$ (characteristic equation) and $x = \lim_{n \rightarrow \infty} \frac{t_{n+1}}{t_n} = 1.38$ (limiting ratio). See also Wilson's {\it Meru 4}.}
\label{table-a1b3gna}
\end{table}

\begin{table}[]
\centering
\begin{tabular}{|l|l|l|l|l|}
\hline
$n$ & creature                       & baby $(b_n)$      & \multicolumn{2}{l|}{total $(t_n)$} \\ \hline
0   & 1                              & 0          & 0                     & 0  \\ \hline
\hline
1   & 1                              & $\alpha$   & $\alpha$              & 2  \\ \hline
2   & 1                              & $\alpha$   & $2\alpha$             & 4  \\ \hline
3   & 1                              & $\alpha$   & $3\alpha$             & 6  \\ \hline
4   & $(1+\alpha)-1$                 & $\alpha^2$ & $\alpha^2+2\alpha$    & 8  \\ \hline
5   & $(\alpha+\alpha)-\alpha$       & $\alpha^2$ & $2\alpha^2+\alpha$    & 10 \\ \hline
6   & $(\alpha+\alpha^2)-\alpha$     & $\alpha^2$ & $3\alpha^2$           & 12 \\ \hline
7   & $(\alpha^2+\alpha^2)-\alpha^2$ & $\alpha^3$ & $\alpha^3+2\alpha^2$  & 16 \\ \hline
8   & $(\alpha^2+\alpha^2)-\alpha^2$ & $\alpha^3$ & $2\alpha^3+\alpha^2$  & 20 \\ \hline
9   & $(\alpha^2+\alpha^3)-\alpha^2$ & $\alpha^3$ & $3\alpha^3$           & 24 \\ \hline
10  & $(\alpha^3+\alpha^3)-\alpha^3$ & $\alpha^4$ & $\alpha^4+2\alpha^3$  & 32 \\ \hline
11  & $(\alpha^3+\alpha^3)-\alpha^3$ & $\alpha^4$ & $2\alpha^4+\alpha^3$  & 40 \\ \hline
12  & $(\alpha^3+\alpha^3)-\alpha^3$ & $\alpha^4$ & $3\alpha^4$           & 48 \\ \hline
\end{tabular}%
\caption{$\alpha =2$, $\beta =2$ and $\gamma =3$. Here $t_{n \ge 2} = b_n + b_{n-1} + b_{n-2}$ (sum of $\beta + 1$ terms). If $n = a \gamma + b$ then $t_n =b \alpha$ (when $a=0$) and $t_n =[b (\alpha - 1) + \gamma] \alpha^a$ (when $a \ge 1$). Also note that $t_{n \le 2} (\alpha =1) = b$ and $t_{n \ge 3} (\alpha =1) = \gamma = 3$.}
\label{table-a2b2g3}
\end{table}

\begin{table}[]
\centering
\begin{tabular}{|l|l|l|l|l|}
\hline
$n$ & creature                        & baby  $(b_n)$      & \multicolumn{2}{l|}{total $(t_n)$} \\ \hline
0   & 1                               & 0           & 0                     & 0  \\ \hline
\hline
1   & 1                               & $\alpha$    & $\alpha$              & 2  \\ \hline
2   & 1                               & $\alpha$    & $2\alpha$             & 4  \\ \hline
3   & 0                               & 0           & $2\alpha$             & 4  \\ \hline
4   & 0                               & 0           & $2\alpha$             & 4  \\ \hline
5   & $\alpha$                        & $\alpha^2$  & $\alpha^2+\alpha$     & 6  \\ \hline
6   & $\alpha+\alpha$                 & $2\alpha^2$ & $3\alpha^2$           & 12 \\ \hline
7   & $(2\alpha +0)-\alpha$           & $\alpha^2$  & $4\alpha^2$           & 16 \\ \hline
8   & $(\alpha+0)-\alpha$             & 0           & $4\alpha^2$           & 16 \\ \hline
9   & $(0+\alpha^2)-0$                & $\alpha^3$  & $\alpha^3+3\alpha^2$            & 20  \\ \hline
10  & $(\alpha^2+2\alpha^2)-0$        & $3\alpha^3$ & $4\alpha^3+\alpha^2$  & 36 \\ \hline
11  & $(3\alpha^2+\alpha^2)-\alpha^2$ & $3\alpha^3$ & $7\alpha^3$           & 56 \\ \hline
12  & $(3\alpha^2+0)-2\alpha^2$       & $\alpha^3$  & $8\alpha^3$           & 64 \\ \hline
13  & $(\alpha^2+\alpha^3)-\alpha^2$  & $\alpha^4$  & $\alpha^4+7\alpha^3$            & 72 \\ \hline
14  & $(\alpha^3+3\alpha^3)-0$      &  $4\alpha^4$     &  $5\alpha^4+4\alpha^3$       & 112  \\ \hline
15  & $(4\alpha^3+3\alpha^3)-\alpha^3$       & $6\alpha^4$        & $11\alpha^4+\alpha^3$                      & 184  \\ \hline
\end{tabular}%
\caption{$\alpha =2$, $\beta =3$ and $\gamma =2$. Here $t_{n \ge 3} = b_n + b_{n-1} + b_{n-2} + b_{n-3}$ (sum of $\beta + 1$ terms).}
\label{table-a2b3g2}
\end{table}

\section{Sums of \(p\)-sequences}
In this section, we obtain closed expressions for odd and even sums, sum of the first \(n\) numbers, and the sum of squares of the first \(n\) numbers of the {\it exponent} \(p\)-sequences. We state the results without giving any proof. We invite the enthusiastic readers to verify (using Tables \ref{table-p2}, \ref{table-p3}, \ref{table-p4} and \ref{table-p5}) and prove them, and also obtain similar expressions for other \(p\)-sequences.

\subsection{Sum of first \(n\) numbers}
Consider the general \(1\)-sequence 
$S_G(1) = \{t_0 = s_0, t_1 = s_0 + a, t_n = t_{n-1} + a = s_0 + na \}$. The sum of first $n \ge 0$ terms of this sequence is 
\be \sum_{k=0}^n t_k(G1) = \frac{(n+1)(t_n + t_0)}{2} =\frac{(n+1)(2s_0 + na)}{2}. \ee

Recall that the {\it exponent} \(p\)-sequence $S_X(p)$ whose seeds are $(0,1,\cdots,p-1)$ is given as 
\be S_X(p) \equiv \{(s_k = k,~ 0\le k \le p-1),~t_n(p)\}, \ee
where $t_n(p) := \sum_{k=n-p}^{n-1} t_k.$ Here we obtain the sum of first \(n\) numbers of \(p\)-sequence $S(Xp)$ for different values of $p$.

\bea
\sum_{k=1}^{n} t_k(X2) &=& t_{n+2} - t_2. \label{eq-sum2x} \\
2 \sum_{k=1}^{n} t_k(X3) &=& t_{n+2} +t_n -[t_3 - t_1]. \label{eq-sum3x} \\
3 \sum_{k=1}^{n} t_k(X4) &=& t_{n+2} +2t_n + t_{n-1} -[t_3 - t_1]. \label{eq-sum4x}
\eea
and

\bea
&&(p-1) \sum_{k=1}^{n} t_k(Xp) \nonumber \\ 
&=& \sum_{k=0}^{p-1} (p-k) t_{n-k} - \big(t_p - \sum_{k=1}^{p-2} (p-1-k) t_{k} \big), \\
&=& t_{n+1} + \sum_{k=0}^{p-2} (p-1-k) t_{n-k} - \big(t_p - \sum_{k=1}^{p-2} (p-1-k) t_{k} \big), \\
&=& t_{n+2} + \sum_{k=0}^{p-3} (p-2-k) t_{n-k} - \big(t_p - \sum_{k=1}^{p-2} (p-1-k) t_{k} \big), \\
&=& t_{n+p} + \sum_{k=1}^{p-2} (p-1-k) t_{n+k} - \big(t_p - \sum_{k=1}^{p-2} (p-1-k) t_{k} \big).
\eea

\subsection{Odd and even sums}
\begin{itemize}

\item Sum of the first \(n\) odd numbers of\\ 
(i) $S_X(2)$ is $\sum_{i=1}^n t_{2i-1} = t_{2n}$, \\~\\
(ii) $S_X(3)$ is $\sum_{i=1}^n t_{2i-1} = \frac{t_{2n} + t_{2n-1} - 1}{2}$, \\~\\
(iii) $S_X(4)$ is $\sum_{i=1}^n t_{2i-1} = \frac{t_{2n+1} + t_{2n-2} - 2}{3}$, and \\~\\
(iv) $S_X(5)$ is $\sum_{i=1}^n t_{2i-1} = \frac{t_{2n+1} + t_{2n-1} + t_{2n-2} + t_{2n-3}}{4}$.\\

\item Sum of the first \(n\) even numbers of\\ 
(i) $S_X(2)$ is $\sum_{i=1}^n t_{2i} = t_{2n+1} - 1$, \\~\\
(ii) $S_X(3)$ is $\sum_{i=1}^n t_{2i} = \frac{t_{2n+1} + t_{2n} - 1}{2}$, \\~\\
(iii) $S_X(4)$ is $\sum_{i=1}^n t_{2i} = \frac{3t_{2n} + 2t_{2n-1}}{3}$, and \\~\\
(iv) $S_X(5)$ is $\sum_{i=1}^n t_{2i} = \frac{4t_{2n} + 2t_{2n-1} + t_{2n-2}}{4}$.


\end{itemize}

\section{Sum of squares of first \(n\) numbers}
Following, we obtain the sum of squares of first $n$ terms of \(p\)-sequence.
For $p=1$ we consider the general sequence $S_G(1)$, and for $p=2,~3,~4$ we consider the {\it exponent} sequence $S_X(p)$.

\subsection{$p=1$ (general sequence)}
For the general \(1\)-sequence 
$S_G(1) = \{t_0 = s_0, t_1 = s_0 + a, t_n = t_{n-1} + a = s_0 + na \}$, the sum of squares of first $n \ge 0$ terms is 
\be \sum_{k=0}^n t_k^2(G1) =(n+1) \big[s_0^2 + ns_0a + \frac{n(2n+1)}{6}a^2 \big]. \label{eq-sum2xsq} \ee

\subsection{$p=2$ (exponent sequence)}
\bea
\sum_{k=1}^n t_k^2(X2) &=& t_n(X2) t_{n+1}(X2). 
\label{eq-sum2xsq}
\eea

\subsection{$p=3$ (exponent sequence)}
\bea
\sum_{k=1}^n t_k^2(X3) &=& t_n(X3) t_{n+1}(X3) - t_{n+1}^2[S_2(3)]. 
\label{eq-sum3xsq}
\eea

\subsection{$p=4$ (exponent sequence)}
\be \sum_{k=1}^n t_k^2(X4) = t_n(X4) t_{n+1}(X4) - \big[(t_{n+2}[S_3(4)] + t_{n-1}[S_3(4)])^2 \pm \delta_n \big], \label{eq-sum4xsq} \ee
where $0 \le \delta_n < t_{n-1}(X4)$. Note that here, unlike $p=2~\&~ 3$ cases, the formula is not exact.

\begin{table}[]
\centering
\begin{tabular}{|l|l|l|}
\hline
$n$ & $\sum_{k=1}^n t_k^2(X2)$                    & $t_n(X2)t_{n+1}(X2)$ \\ \hline
1   & $1^2$                            & $1 \times 1$    \\ \hline
2   & $1 + 1^2 = 2$ & $1 \times 2$    \\ \hline
3   & $2 + 2^2 = 6$                    & $2 \times 3$    \\ \hline
4   & $6 + 3^2 = 15$                   & $3 \times 5$    \\ \hline
5   & $15 + 5^2 = 40$                  & $5 \times 8$    \\ \hline
6   & $40 + 8^2 = 104$                 & $8 \times 13$   \\ \hline
7   & $104 + 13^2 = 273$               & $13 \times 21$  \\ \hline
8   & $273 + 21^2 = 714$               & $21 \times 34$  \\ \hline
9   & $714 + 34^2 = 1870$              & $34 \times 55$  \\ \hline
\end{tabular}%
\caption{Sum of squares of numbers of \(2\)-sequence $S_X(2)$ 
(see Table \ref{table-p2}).
}
\label{table-sum2xsq}
\end{table}

\begin{table}[]
\centering
\begin{tabular}{|l|l|l|}
\hline
$n$ & $\sum_{k=1}^n t_k^2(X3)$                    & $t_n(X3) t_{n+1}(X3) - t_{n+1}^2[S_2(3)]$      \\ \hline
1   & $1^2$                            & $1 \times 2 - 1^2$      \\ \hline
2   & $1 + 2^2 = 5$ & $2 \times 3 -  1^2$     \\ \hline
3   & $5 + 3^2 = 14$                   & $3 \times 6 - 2^2$      \\ \hline
4   & $14 + 6^2 = 50$                  & $6 \times 11 - 4^2$    \\ \hline
5   & $50 + 11^2 = 171$                & $11 \times 20 - 7^2$    \\ \hline
6   & $171 + 20^2 = 571$               & $20 \times 37 - 13^2$   \\ \hline
7   & $571 + 37^2 = 1940$              & $37 \times 68 - 24^2$   \\ \hline
8   & $1940 + 68^2 = 6564$             & $68 \times 125 - 44^2$  \\ \hline
9   & $6564 + 125^2 = 22189$           & $125 \times 230 - 81^2$ \\ \hline
\end{tabular}%
\caption{Sum of squares of numbers of \(3\)-sequence $S_X(3)$ 
(see Table \ref{table-p3}).
}
\label{table-sum3xsq}
\end{table}

\begin{table}[]
\centering
\begin{tabular}{|l|l|l|}
\hline
$n$ & $\sum_{k=1}^n t_k^2(X4)$                    & $t_n(X4) t_{n+1}(X4) - \big[(t_{n+2}[S_3(4)] + t_{n-1}[S_3(4)])^2 \pm \delta_n \big]$      \\ \hline
1   & $1^2$                            & $1 \times 2 - [(1+0)^2 + 0]$      \\ \hline
2   & $1 + 2^2 = 5$ 					& $2 \times 3 -  [(1+0)^2 + 0]$     \\ \hline
3   & $5 + 3^2 = 14$                   & $3 \times 6 - [(2+0)^2 + 0]$      \\ \hline
4   & $14 + 6^2 = 50$                  & $6 \times 12 - [(4+1)^2 - 3]$    \\ \hline
5   & $50 + 12^2 = 194$                & $12 \times 23 - [(8+1)^2 + 1]$    \\ \hline
6   & $194 + 23^2 = 723$               & $23 \times 44 - [(15+2)^2 + 0]$   \\ \hline
7   & $723 + 44^2 = 2659$              & $44 \times 85 - [(29+4)^2 - 8]$   \\ \hline
8   & $2659 + 85^2 = 9884$             & $85 \times 164 - [(56+8)^2 - 40]$  \\ \hline
9   & $9884 + 164^2 = 36780$           & $164 \times 316 - [(108+15)^2 + 67]$ \\ \hline
\end{tabular}%
\caption{Sum of squares of numbers of \(4\)-sequence $S_X(4)$ (see Table \ref{table-p4}).
}
\label{table-sum4xsq}
\end{table}

\subsection{More sums}
We further find that
\begin{itemize}

\item For \(2\)-sequence $S_X(2)$,\\ 
(i) $\sum_{i=1}^n t_i t_{i+1} = \frac12 (t_{n+2}^2 - t_nt_{n+1} - t_2^2)$, and\\
(ii) $\sum_{i=1}^n t_i t_{i+1} = \sum_{j=1}^n (n+1-j) t_j^2$.\\

\item For \(3\)-sequence $S_X(3)$,\\ 
$\sum_{i=1}^{n-1} (t_i t_{i+1} + t_{i+1} t_{i+2} + t_{i+2} t_{i}) \\ = \frac12 \big[t_{n+2}(t_{n+2} -1) + t_n(t_n - 1) - (t_3 - t_1)(t_3 + t_1 -1) \big]$.\\

\end{itemize}

\section{\(p\)-golden ratio}
In this section, we present golden ratio ($\Phi_p$) and golden angle ($\theta_g(p)$) associated with \(p\)-sequences, and consider other generalizations of golden ratio.

The {\it golden ratio} is one of the most famous numbers. Given $a$ and $b(<a)$ two positive numbers, the golden ratio is defined as 
\be \label{eq-2ratio} \frac{a}{b} = \frac{a+b}{a}. \ee 
Taking $\frac{a}{b} = \Phi$, 
Eq. (\ref{eq-2ratio}) reduces to the quadratic equation $\Phi^2 = \Phi +1$ whose positive solution is $\Phi = \frac{\sqrt{5}+1}{2} = 1.61803$. This value corresponds to the limiting ratio value of the Fibonacci sequence. 

For the golden ratios associated with \(p\)-sequences, we first ask a couple of questions: {\it (i) Does there exist a ratio, like golden ratio [Eq. (\ref{eq-2ratio})], for given \(p \ge 3\) positive real numbers? (ii) What is the value of this ratio? Is this value unique? 
(ii) Is this value of ratio equal to the limiting ratio value of \(p\)-sequences?} Surprisingly enough, the answer is in affirmative. 

\begin{figure}%
\centering
\includegraphics[width = 3in]{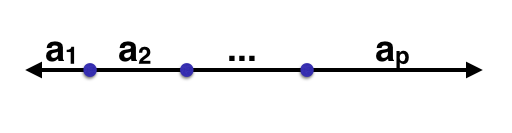}%
\caption{Division of a line into \(p\) segments.}
\label{fig-goldenratiop}
\end{figure}

Suppose $a_1<a_2<\cdots<a_p$ are \(p \ge 2\) positive real numbers (see Fig. \ref{fig-goldenratiop}). We define the \(p\)-golden ratio as
\footnote{We will see later that actually $\frac{t_{n+1}}{t_n}$ is the golden ratio for large \(n\). The relation of limiting ratio value of \(p\)-sequence with the Euclid's problem, Eq. (\ref{eq-p-ratio}), is accidental.} 
\be \label{eq-p-ratio}
\frac{a_2}{a_1} = \frac{a_3}{a_2} = \cdots = \frac{\sum_{k=1}^p a_k}{a_p} (= \Phi_p). \ee
Note that Eq. (\ref{eq-2ratio}) is a special case of Eq. (\ref{eq-p-ratio}) for $p=2$.

\subsection{Characteristic equation for $\Phi_p$}
We find that from Eq. (\ref{eq-p-ratio}) follows naturally the $p$-degree algebraic equation whose positive solution gives the value of $\Phi_p$:  
\be \label{eq-godeqn}
X_p(x) \equiv x^p - \sum_{k=0}^{p-1} x^k = 0. \ee
We call this {\it golden equation}.
Note that $X_p(0) = -1$ for all \(p\) and $X_p(1) = -(p-1)$.
This equation has been obtained recently in an interesting physical problem concerning center of masses in two and higher dimensions \cite{dutta2020}.

\subsection{Less radical characteristic equations}
For fixed \(p\) and positive integers $\{k_i\}$, one can choose recurrence relations with $m$ terms ($2 \le m,k_i < p$) to obtain the following less radical characteristic equations, 
\bea 
x^p &=& 1 + x^{k_1}, \nonumber \\
x^p &=& 1 + x^{k_1} + x^{k_2}, \nonumber \\
x^p &=& 1 + x^{k_1} + x^{k_2} + x^{k_3}, \nonumber
\eea
and so on, each with its own convergence. Wilson's {\it Meru 1} through {\it Meru 9} are particular examples of the above \(2\)-term characteristic equation for $p=2,3,4,5$.

\subsection{Roots of the golden equation}
\begin{figure}%
\centering
\includegraphics[width = 4.5in]{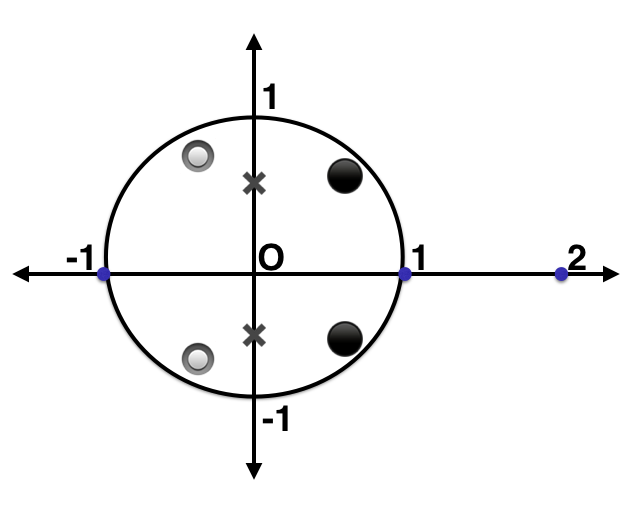}%
\caption{Schematic representation of roots of the {\it golden equation}: $x^p = \sum_{k=0}^{p-1} x^k$.}
\label{fig-goldenroots1}
\end{figure}

\begin{figure}
\subfloat[$p=1$]{\includegraphics[width = 2.7in]{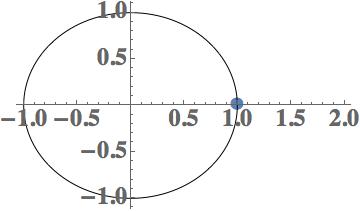}}~~~
\subfloat[$p=2$]{\includegraphics[width = 2.7in]{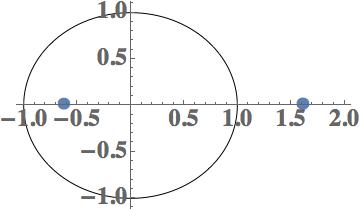}}\\
\subfloat[$p=3$]{\includegraphics[width = 2.7in]{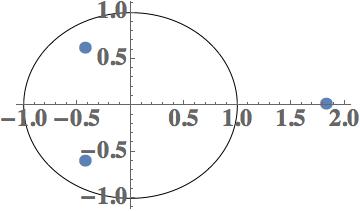}}~~~
\subfloat[$p=5$]{\includegraphics[width = 2.7in]{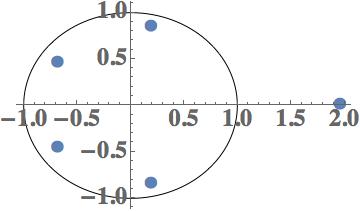}} \\
\subfloat[$p=10$]{\includegraphics[width = 2.7in]{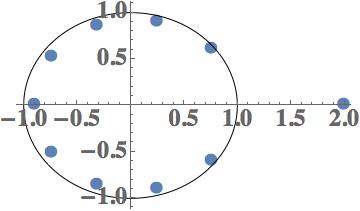}}~~~ 
\subfloat[$p=15$]{\includegraphics[width = 2.7in]{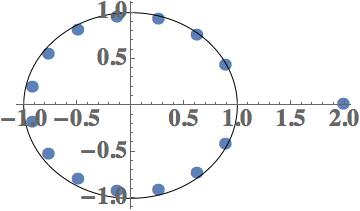}}\\
\subfloat[$p=18$]{\includegraphics[width = 2.7in]{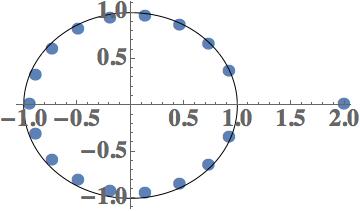}}~~~
\subfloat[$p=24$]{\includegraphics[width = 2.7in]{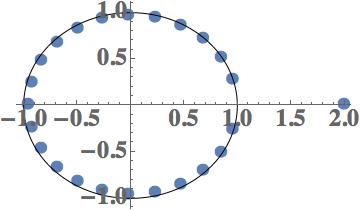}} 
\caption{Roots of the {\it golden equation} for different values of \(p\).}
\label{fig-goldenroots2}
\end{figure}

Here we look at the nature of roots of Eq. (\ref{eq-godeqn}). Roots can be positive, negative and complex. Complex roots obviously occur in pairs and lie within a unit circle and approaches towards the boundary of the circle with increasing \(p\). The only negative root approaches $-1$ for large \(p\). The only positive root lies between \(1\) and \(2\), and tends to $2$ for large \(p\). See Figs. \ref{fig-goldenroots1} and \ref{fig-goldenroots2}.

\subsection{Near cousins of the golden equation}
\begin{figure}%
\subfloat[$x^p = x^{p-1} - x^{p-2} + \cdots \pm 1$]{\includegraphics[width = 3.0in]{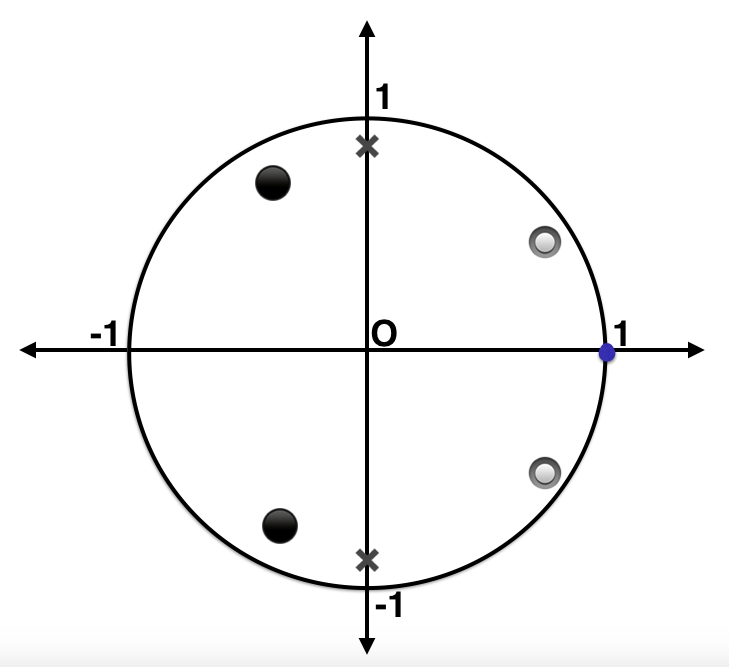}}%
\subfloat[$x^p = - \sum_{k=0}^{p-1} x^k$]{\includegraphics[width = 3.0in]{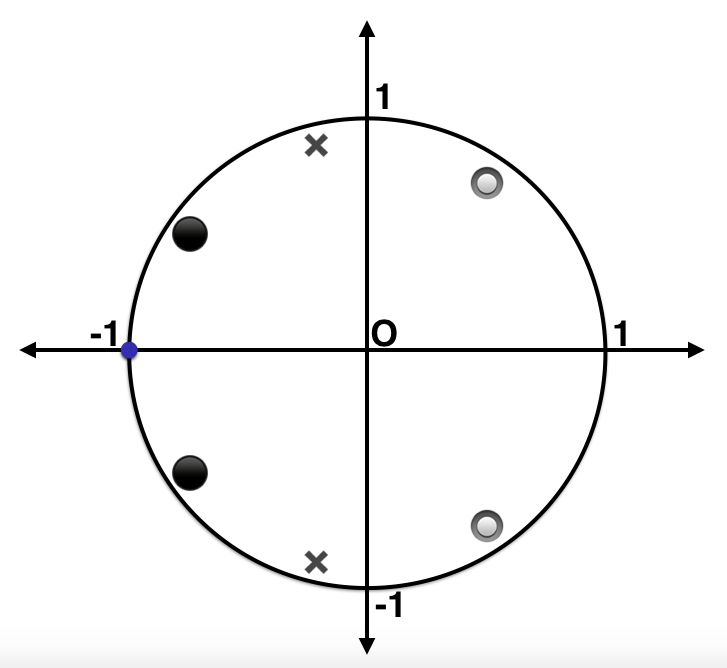}}%
\caption{Schematic representation of roots of the near cousins of the golden equation. 
}
\label{fig-cousins}
\end{figure}

The equations $(a)~x^p = x^{p-1} - x^{p-2} + \cdots \pm 1$ and $(b)~x^p = - \sum_{k=0}^{p-1} x^k$ are two immediate near cousins of the golden equation $x^p = \sum_{k=0}^{p-1} x^k$. Roots of both $(a)~\text{and}~(b)$ are complex and real. Complex roots obviously occur in pairs and lie within a unit circle and approaches towards the boundary of the circle with increasing \(p\). The only real positive root of $(a)$ is $+1$, and the only real negative root of $(b)$ is $-1$. 
See Fig. \ref{fig-cousins}.

\section{Recursion relation for $\Phi_p$}
Because $\Phi_p$ is a solution of Eq. (\ref{eq-godeqn}), we have 
\bea 
\Phi_p^p &=& \Phi_p^{p-1} + \Phi_p^{p-2} + \cdots + \Phi_p + 1 = \sum_{k=0}^{p-1} \Phi_p^k, 
\label{eq-basic} \\
\Phi_p^{p+1} &=& \Phi_p^{p} + \Phi_p^{p-1} + \cdots + \Phi_p^2 + \Phi_p, \nonumber \\
&=& 2\Phi_p^{p} - 1. \label{eq-basic0}
\eea
Eq. (\ref{eq-basic}) can be equivalently rewritten as
\bea 
\Phi_p &=& 1 + \frac{\sum_{k=0}^{p-2} \Phi_p^k}{\Phi_p^{p-1}}, \label{eq-basic1} \\
&=& 1 + \frac{1}{\Phi_p-1 + \frac{1}{\sum_{k=0}^{p-2} \Phi_p^k}}. \label{eq-basic2}
\eea

Also, Eq. (\ref{eq-basic}) implies a recursion relation
\be \label{eq-recursion}
\Phi_p^n = \Phi_p^{n-1} + \Phi_p^{n-2} + \cdots + \Phi_p^{n-p} = \sum_{k=n-p}^{n-1} \Phi_p^k.
\ee

\section{$\Phi_1$ of \(1\)-sequence}
We have seen above that Eq. (\ref{eq-basic}) is the basic equation for $\Phi_{p\ge2}$.
If we consider this sacred golden equation for $p=1$, we have
\be \label{eq-p1} \Phi_1 = \Phi_1^0 = 1. \ee 
We remark that $\Phi_1$ is related with the limiting ratio value of \(1\)-sequences. We construct a \(1\)-sequence by choosing a seed $s_0 \ge 0$ and a constant $a \ge 0$ such that $t_0 = s_0$, and for $n \ge 1$
\be t_n = t_{n-1} +a = s_0 + na. \ee 
The limiting ratio value for this \(1\)-sequence is then 
\be \label{eq-1seq-ratio-limit}
\lim_{n \rightarrow \infty} \frac{t_{n+1}(1)}{t_{n}(1)} = \lim_{n \rightarrow \infty} \frac{s_0 + (n+1)a}{s_0 + na} = \lim_{n \rightarrow \infty} \big(1+ \frac{1}{n + \frac{s_0}{a}} \big) = 1 = \Phi_1. \ee

We note that Eq. (\ref{eq-p1}) provides the {\it lower limit} on $\Phi_p$'s. That is,
\be \label{eq-lowerbound} \Phi_{p \ge 1} \ge 1. \ee

\begin{table}[htb]
\centering
\begin{tabular}{|l|l|l|l|}
\hline
$p$ & $\Phi_p$ & $\theta = \arcsin[\frac{\Phi_p-1}{2}]$ & $\theta = \arcsin[\Phi_p/2]$ \\ \hline
\hline
1   & 1.0 & 0.0 & 30.0                              \\ \hline
\hline
2   & 1.61803   & 18.0    & 54.0                   \\ \hline
3   & 1.83929   & 24.8122    & 66.8742           \\ \hline
4   & 1.92756   & 27.6313    & 74.5321     \\ \hline
5   & 1.96595   & 28.8799    & 79.4124       \\ \hline
6   & 1.98358   & 29.4583    & 82.6531          \\ \hline
7   & 1.99196   & 29.7344    & 84.8608        \\ \hline
8   & 1.99603   & 29.8688    & 86.3893              \\ \hline
9   & 1.99803   & 29.9349  & 87.4567       \\ \hline
10  & 1.99902  & 29.9676    & 88.2063         \\ \hline
11  & 1.99951   & 29.9838    & 88.7317         \\ \hline
12  & 1.99976   & 29.9921    & 89.1124               \\ \hline
13  & 1.99988   & 29.996    & 89.3724                \\ \hline
14  & 1.99994   & 29.998    & 89.5562             \\ \hline
15  & 1.99997   & 29.999    & 89.6862                \\ \hline
16  & 1.99998   & 29.9993  & 89.7438                      \\ \hline
17  & 1.99999   & 29.9997  & 89.8188                  \\ \hline
\hline
18  & 2.0       & 30.0   & 90.0                            \\ \hline
19  & 2.0       & 30.0   & 90.0                                \\ \hline
20  & 2.0       & 30.0   & 90.0                               \\ \hline
\end{tabular}%
\caption{Values of $\Phi_p$, and the trigonometric angles such that $\Phi_p = 1 + 2\sin \theta = 2\sin \theta$.}
\label{table-phi}
\end{table}

\begin{figure}%
\centering
\resizebox{0.6\columnwidth}{!}{
\includegraphics{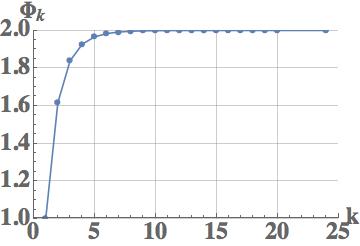}%
}%
\caption{Plot of \(p\)-golden ratios.}
\label{fig-goldenratio plot}
\end{figure}

\newpage
\section{$\Phi_p^n$ as a polynomial of degree $p-1$}  
We have seen earlier the following relations:
$$\Phi_p^p = \Phi_p^{p-1} + \Phi_p^{p-2} + \cdots + \Phi_p + 1,$$ and
$$\Phi_p^n = \Phi_p^{n-1} + \Phi_p^{n-2} + \cdots + \Phi_p^{n-p}.$$ 
Here the question we want to address is: {\it is it possible to reduce $\Phi_p^n$ to a polynomial of degree $p-1$?} Put differently, can we express $\Phi_p^n$ in terms of $\{\Phi_p^k\}_{k=0}^{p-1}$? 
It is very illuminating to see that it is possible to express $\Phi_p^n~(n \ge 0)$ in terms of $\{\Phi_p^k\}_{k=0}^{p-1}$ as follows:
\bea 
\Phi_p^n  & =& t_n[S_{p-1}(p)] \Phi_p^{p-1} + \cdots + t_n[S_{1}(p)] \Phi_p +t_n[S_{0}(p)]
 \nonumber \\
&=& \sum_{k=0}^{p-1} t_n[S_{k}(p)] \Phi_p^{k},
\label{eq-phin}
\eea
where $t_n[S_{k}(p)] = \sum_{j=n-p}^{n-1} t_j[S_{k}(p)]$.
Eq. (\ref{eq-phin}) can be easily verified from Tables \ref{table-p2}, \ref{table-p3}, \ref{table-p4} and \ref{table-p5}.\\ 

In particular, for $p=2~\text{and}~3$, the explicit expressiosn are

\be\label{eq:phi2}
\Phi_2^n = \left\{
  \begin{array}{ll}
   t_n[S_1(2)]\Phi_2 + t_n[S_0(2)] \quad &(n \ge 0),\\
   t_{n+1}[S_0(2)]\Phi_2 + t_n[S_0(2)] \quad &(n \ge 2),\\
   t_{n}[S_X(2)]\Phi_2 + t_{n-1}[S_X(2)] \quad &(n \ge 2),
  \end{array}
\right.
\ee

and

\be\label{eq:phi3}
\Phi_3^n = \left\{
  \begin{array}{lr}
   t_n[S_2(3)]\Phi_3^2 + t_n[S_1(3)]\Phi_3 + t_n[S_0(3)] \quad &(n \ge 0),\\
   t_{n-3}[S_S(3)]\Phi_3^2 + t_{n-2}[S_X(3)]\Phi_3 + t_{n-4}[S_S(3)] \quad &(n \ge 4).
  \end{array}
\right.
\ee

\section{Applications of \(p\)-golden ratios}
We have seen earlier that the golden ratio and the related Fibonacci sequence are present in abundance in our everyday life. We also learnt the skeptical view on this, and that not all objects exhibit the golden ratio in the sense that convergent limits do not settle down to the numerical value $1.618$. This is now evident with the introduction of \(p\)-sequences and the associated \(p\)-golden ratios why it is not the case. In fact, $\Phi_2 = 1.618$ is only one member of several families of golden ratios (such as those of Stakhov, Spinadel, Krcadinac, etc. including the present work). Therefore, it is natural to expect that $\Phi_{p>2}$ will have many interesting applications as well. 

\section{Golden geometry}
\subsection{Golden angles}
\begin{figure}%
\subfloat[]{\includegraphics[width = 2.4in]{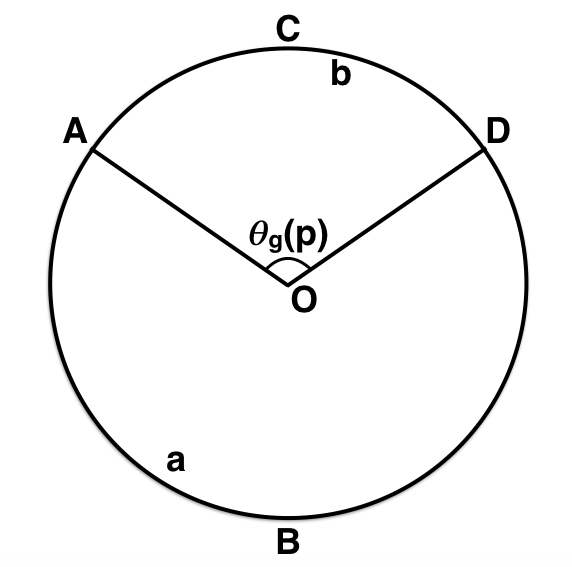}}~~~
\subfloat[]{\includegraphics[width = 3.0in]{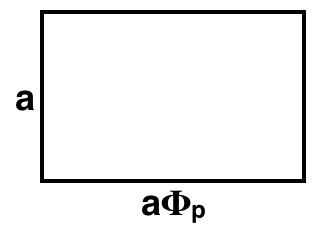}}
\caption{The golden angle $\theta_g$ is determined by using $a/b = \Phi_p$.}
\label{fig-goldenang}
\end{figure}

The golden angle is defined as the acute angle $\theta_g$ that divides the circumference of a circle into two arcs \(ABD\) and \(ACD\) with lengths in the golden ratio. See Fig. \ref{fig-goldenang}(a). The golden ratio here satisfies $\Phi_p = \frac{a}{b}$. We then determine the golden angle by $\frac{\theta_g(p)}{2\pi} = \frac{b}{a+b} = \frac{1}{1+\frac{a}{b}} = \frac{1}{1+\Phi_p}$. Hence,
\be \theta_g(p) = \frac{2\pi}{1+\Phi_p}. \label{eq-goldenang} \ee
From Table \ref{table-goldenang} we see that $\frac{2\pi}{3} \le \theta_g(p) \le \pi$.

\subsection{Golden shapes}
We can construct geometrical objects such as polygons (rectangle, pentagon, etc.) and spirals which have properties characterizing the golden \(p\)-ratio or certain \(p\)-sequences. Note that {\it a square is a golden rectangle with golden ratio $\Phi_1=1$}. 



\begin{table}[]
\centering
\begin{tabular}{|l|l|l|l|l|}
\hline
$p$       & 1       & 2       & 3       & 18       \\ \hline
$\Phi_p$       & 1.0       & 1.61803       & 1.83929       & 2.0       \\ \hline
$\theta_g(p)$ & $180^{\circ}$ & $137.5^{\circ}$ & $126.8^{\circ}$ & $120^{\circ}$ \\ \hline
\end{tabular}%
\caption{The golden angles for \(p\)-sequences, $p=1,~2,~3,~18$.}
\label{table-goldenang}
\end{table}

\section{Further generalizations of golden ratio}
The trouble with the notion of golden ratio is that it can be extended in many ways such that the original golden ratio $\Phi_2$ is a particular case.
In an earlier section, we have seen that the recurrence relation
$t_n(p) = \sum_{k=1}^p t_{n-k}(p)$
and the golden ratio 
$\frac{a_2}{a_1} = \frac{a_3}{a_2} = \cdots = \frac{\sum_{k=1}^p a_k}{a_p}$
correspond to the characteristic equation
$x^p = \sum_{k=0}^{p-1} x^k$.
A straight-forward generalization of these yield
\bea
t_n(p) &=& \sum_{k=1}^p c_k t_{n-k}(p), \label{eq-grr1}\\
\frac{a_2}{a_1} &=& \frac{a_3}{a_2} = \cdots = \frac{\sum_{k=1}^p c_k a_k}{a_p}, \label{eq-ggr1}\\
x^p &=& \sum_{k=0}^{p-1} c_k x^k. \label{eq-gce1}
\eea

That is, for a sequence of numbers whose terms are given by the (weighted) sum of its {\it consecutive \(p\)-previous terms}, the characteristic polynomial equation can be obtained by using the golden ratio. However, how do we obtain the characteristic polynomial equation for an arbitrary recurrence relation
\footnote{Wilson's {\it Meru 1} through {\it Meru 9} \cite{wilson1993, wilson2001} with their limiting ratios are particular examples of Eq. (\ref{eq-grr2}). 
}, 
\be \label{eq-grr2}
t_n = c_1 t_{n-m_1} + c_2 t_{n-m_2} + \cdots + c_p t_{n-m_p}, \ee
otherwise? 
In this case also, we can project a ratio like the golden one, 
Eq. (\ref{eq-ggr1}), as given below
\be x = \frac{t_{n-m+1}}{t_{n-m}} = \frac{t_{n-m+2}}{t_{n-m+1}} = \cdots = \frac{t_{n}}{t_{n-1}}, \ee
where $m = \text{max} \{m_1,~m_2,~\cdots,~m_p\}$ so that 
\bea 
t_{n-m_k} &=& x^{m-m_k} t_{n-m},~~(1 \le k \le p) \nonumber \\
t_{n-1} &=& x^{m-1} t_{n-m}.
\eea
Then, the characteristic polynomial equation is
\footnote{Proof of Eq. (\ref{eq-gce2}).
\bea
x t_{n-1} = t_n &=& c_1 t_{n-m_1} + c_2 t_{n-m_2} + \cdots + c_p t_{n-m_p}, \nonumber \\
\Rightarrow x (x^{m-1} t_{n-m}) &=& (c_1 x^{m-m_1} + c_2 x^{m-m_2} + \cdots + c_p x^{m-m_p})
t_{n-m},
\nonumber \\
\Rightarrow x^{m} &=& c_1 x^{m-m_1} + c_2 x^{m-m_2} + \cdots + c_p x^{m-m_p}. \nonumber
\eea
}
\be \label{eq-gce2} 
x^m = c_1 x^{m-m_1} + c_2 x^{m-m_2} + \cdots + c_p x^{m-m_p}. \ee  
We state a proposition below which gives us a straightforward general rule to obtain the characteristic polynomial equation for an arbitrary recurrence relation.\\
{\it Proposition. The polynomial equation characteristic to a given recurrence relation is obtained by requiring $x^{u-v} := \lim_{n \rightarrow \infty} \frac{t_{n+u}}{t_{n+v}}$, where \(u\) and \(v\) are integers. The characteristic equation is the minimal polynomial which gives the value of the limiting ratio of the sequence, and from which all its algebraic properties follow. For the generalized recurrence relation, $t_n = c_1 t_{n-m_1} + c_2 t_{n-m_2} + \cdots + c_p t_{n-m_p}$, the characteristic polynomial equation is given by $x^m = c_1 x^{m-m_1} + c_2 x^{m-m_2} + \cdots + c_p x^{m-m_p}$, where $m = \text{max} \{m_1,~m_2,~\cdots,~m_p\}$
\footnote{Another proof of Eq. (\ref{eq-gce2}).
\bea 
x &=& \lim_{n \rightarrow \infty} \frac{t_{n+1}}{t_n}, \nonumber \\
&=& \lim_{n \rightarrow \infty} \frac{c_1 t_{n-(m_1-1)} + c_2 t_{n-(m_2-1)} + \cdots + c_p t_{n-(m_p-1)}}{t_n}, \nonumber \\
&=& c_1 \lim_{n \rightarrow \infty} \frac{t_{n-(m_1-1)}}{t_n} + c_2 \lim_{n \rightarrow \infty} \frac{t_{n-(m_2-1)}}{t_n} + \cdots + c_p \lim_{n \rightarrow \infty} \frac{t_{n-(m_p-1)}}{t_n}, \nonumber \\
&=& c_1 x^{-(m_1-1)} + c_2 x^{-(m_2-1)} + \cdots + c_p x^{-(m_p-1)}, \nonumber \\
&=& \frac{c_1 x^{m-m_1} + c_2 x^{m-m_2} + \cdots + c_p x^{m-m_p}}{x^{m-1}}, \nonumber \\
\Rightarrow x^{m} &=& c_1 x^{m-m_1} + c_2 x^{m-m_2} + \cdots + c_p x^{m-m_p}. \nonumber
\eea
Thus, $\lim_{n \rightarrow \infty} \frac{t_{n+1}}{t_{n}}$ is the golden ratio in general.}
}.\\~\\

Moving a step further, we consider the relation
\be \big(u_1 \frac{a_2}{a_1} \big)^{v_1} = \big(u_2 \frac{a_3}{a_2} \big)^{v_2} = \cdots = \big(u_{p-1} \frac{a_p}{a_{p-1}} \big)^{v_{p-1}} = \big(u_p \frac{\sum_{k=1}^p c_k a_k}{a_p} \big)^{v_p}, \label{eq-gen-golden-ratio} \ee
where $\{(u_i,~v_i)\}$ and $\{c_k\}$ are given. Goal is to find values of the ratios $\{\frac{a_{k+1}}{a_k}\}$ and $\frac{\sum_{k=1}^p c_k a_k}{a_p}$ such that Eq. (\ref{eq-gen-golden-ratio}) holds. Does a solution exist? This problem is rather hard to solve in general.

Next, one can choose any pair of ratios at a time. Say, $\big(u_1 \frac{a_2}{a_1} \big)^{v_1} = \big(u_2 \frac{a_3}{a_2} \big)^{v_2}$. There are two cases here. (i) Assume that $\frac{a_2}{a_1} = x$ and $\frac{a_3}{a_2} = f_{23}(x)$. Then the characteristic equation is $\big(u_1 x \big)^{v_1} = \big(u_2 f_{23}(x) \big)^{v_2}$ and the positive solution is $x= \frac{1}{u_1} \big(u_2 f_{23}(x) \big)^{\frac{v_2}{v_1}}$. (ii) For $\frac{a_3}{a_2} = x$ and $\frac{a_2}{a_1} = f_{12}(x)$, the characteristic equation is $\big(u_1 f_{12}(x) \big)^{v_1} = \big(u_2 x \big)^{v_2}$ and the positive solution is $x= \frac{1}{u_2} \big(u_1 f_{12}(x) \big)^{\frac{v_1}{v_2}}$. Thus, equating two ratios at a time, we will have $2(p-1)!$ characteristic polynomial equations and consequently as many roots of them for given $\{(u_i,~v_i)\}$ and $\{c_k\}$. 
To the best of our knowledge, most generalizations of the Fibonacci sequence and the golden ratio can be seen as special cases of Eqs. (\ref{eq-grr2}), (\ref{eq-gce2}) and (\ref{eq-gen-golden-ratio}).

\section{\(Q\)-matrix and determinantal identity}
We have seen that the Fibonacci sequence can be generalized in many ways such as generalizations of the Euclid's theorem, the recurrence relations, and the characteristic equations. There is yet another way to study and generalize the Fibonacci sequence and derive many interesting properties of these numbers using a matrix representation \cite{miles1960, ivie1972, silvester1979, kalman1982, er1984a, er1984b, lee1997, stakhov1999, lee2003, karaduman2005, stanimirovic2008, chatterjee2018}.
By matrix methods, while Silvester \cite{silvester1979} derived many interesting properties of the Fibonacci numbers, Kalman \cite{kalman1982} generalized Fibonacci numbers.
In this section, we consider the generating \(Q\)-matrix and the {\it determinantal identities} of \(p\)-sequences.\\

For the Fibonacci sequence recurrence relation $f_n = f_{n-1} + f_{n-2}$ given $f_1=1~\text{and}~f_0=0$, using the generating $Q$-matrix
\begin{equation} \label{eq-q2basic}
Q = 
\begin{pmatrix}
1 & 1 \\
1 & 0 
\end{pmatrix}
 =
 \begin{pmatrix}
f_2 & f_1 \\
f_1 & f_0 
\end{pmatrix},
\end{equation}
we have, for $n \ge 1$,
\begin{equation} \label{eq-q2nseq1}
\begin{pmatrix}
f_{n+1} \\
f_n 
\end{pmatrix}
=
\begin{pmatrix}
1 & 1 \\
1 & 0 
\end{pmatrix}^n
\begin{pmatrix}
f_1 \\
f_0 
\end{pmatrix},
\end{equation}
and
\begin{equation} \label{eq-q2nseq2}
Q^n = 
\begin{pmatrix}
f_{n+1} & f_n \\
f_n & f_{n-1} 
\end{pmatrix}.
\end{equation}

From Eq. (\ref{eq-q2nseq2}) follows the Cassini's identity
\begin{equation} \label{eq-2cassini}
(\text{det}~Q)^n = 
\begin{vmatrix}
f_{n+1} & f_n \\
f_n & f_{n-1} 
\end{vmatrix}
=(-1)^n.
\end{equation}
\begin{table}[]
\centering
\begin{tabular}{|l|l|}
\hline
$n$ & $f_{n+1} f_{n-1} - f_n^2 = (-1)^n$  \\ \hline
1   & $1 \times 0 - 1^2 = -1$                                \\ \hline
2   & $2 \times 1 - 1^2 = 1$                     \\ \hline
3   & $3 \times 1 - 2^2 = -1$                      \\ \hline
4   & $5 \times 2 - 3^2 = 1$                      \\ \hline
5   & $8 \times 3 - 5^2 = -1$                     \\ \hline
6   & $13 \times 5 - 8^2 = 1$                    \\ \hline
7   & $21 \times 8 - 13^2 = -1$                 \\ \hline
8   & $34 \times 13 - 21^2 = 1$                  \\ \hline
9   & $55 \times 21 - 34^2 = -1$                 \\ \hline
\end{tabular}%
\caption{Illustration of the Cassini's identity for the Fibonacci sequence $S_1(2)$ (see Table \ref{table-p2}).}
\label{table-2cassini}
\end{table}

In the following, we confine ourselves to establishing a family of {\it determinantal identities} of which the Cassini's identity is a particular case.
We begin with introducing 
a square matrix $Q_p$ of order \(p\) 
\begin{equation} \label{eq-qpbasic}
Q_{p} = 
\begin{pmatrix}
1 & 1 & 1 & \cdots & 1 \\
1 & 0 & 0 & \cdots & 0 \\
0 & 1 & 0 & \cdots & 0 \\
\vdots  & \vdots  & \ddots & \vdots  & \vdots  \\
0 & \cdots & 0 & 1 & 0 
\end{pmatrix},
\end{equation}
where elements of the first row are seeds of the {\it coefficient} \(p\)-sequence $S_{C}(p)$, those of the second row are seeds of \(p\)-sequence $S_{0}(p)$, the third row are seeds of \(p\)-sequence $S_{1}(p)$, and so on, until the last row whose elements are seeds of \(p\)-sequence $S_{p-2}(p)$. 
The determinant of this matrix $Q_p$ is
\be\label{eq:qpdet}
\text{det}~Q_p = \left\{
  \begin{array}{ll}
   -1 \quad &\text{(\(p\) even)},\\
    1 \quad &\text{(\(p\) odd)}.
  \end{array}
\right.
\ee

The matrix $Q_p$ can be seen as a special case of the generating matrix $\tilde{Q}_p$ \cite{kalman1982},
\begin{equation} \label{eq-qpbasic-gen}
\tilde{Q}_{p} = 
\begin{pmatrix}
c_1 & c_2 & c_3 & \cdots & c_p \\
1 & 0 & 0 & \cdots & 0 \\
0 & 1 & 0 & \cdots & 0 \\
\vdots  & \vdots  & \ddots & \vdots  & \vdots  \\
0 & \cdots & 0 & 1 & 0 
\end{pmatrix},
\end{equation}
for the generalized recurrence relation $t_n(p) = \sum_{k=1}^p c_k t_{n-k}(p)$, where $\{c_k\}$s are constants, such that 
\begin{equation} \label{eq-q2nseq1}
\begin{pmatrix}
t_{n+p-1} \\
\vdots \\
t_{n+1} \\
t_n 
\end{pmatrix}
=
\tilde{Q}^n_p
\begin{pmatrix}
t_{p-1} \\
\vdots \\
t_1 \\
t_0 
\end{pmatrix}
~~(n \ge 1),
\end{equation}
and 
\begin{equation}
\text{det}~\tilde{Q}_p = (-1)^{p+1}c_p.
\end{equation}

We can, further, identify the matrix $Q_p$ as
\begin{equation} \label{eq-qpterms}
Q_{p} = 
\begin{pmatrix}
t_{p} & t_{p-1} & \cdots & t_{p-1} & t_{p-1} \\
t_{p-1} & 0 & \cdots & 0 & t_{p-2} \\
t_{p-2} & t_{p-1} & \cdots & 0 & t_{p-3} \\
\vdots  & \vdots  & \ddots & \vdots  & \vdots  \\
t_{1} & 0 & \cdots & t_{p-1} & t_{0} 
\end{pmatrix},
\end{equation}
where $t_k$'s in the first row, the first and the last columns, and the diagonal below the main diagonal are terms of the \(p\)-sequence $S_{p-1}(p)$, and other elements of the matrix are zero. 
Now, making use of Eqs. (\ref{eq-qpbasic}) and (\ref{eq-qpterms}), we can construct the determinantal identities, like the Cassini's identity, for $p \ge 3$. Following, we illustrate the case for $p=3,~4$.

\subsection{Illustration for $p=3$}
Starting with
\begin{equation}
Q_{3} = 
\begin{pmatrix}
1 & 1 & 1 \\
1 & 0 & 0 \\
0 & 1 & 0
\end{pmatrix}
=
\begin{pmatrix}
t_{3} & t_{2} & t_{2} \\
t_{2} & 0 & t_{1} \\
t_{1} & t_{2} & t_{0}
\end{pmatrix},
\end{equation}

we have
\bea
Q^2_{3} &=& 
\begin{pmatrix}
1 & 1 & 1 \\
1 & 0 & 0 \\
0 & 1 & 0
\end{pmatrix}
\begin{pmatrix}
t_{3} & t_{2} & t_{2} \\
t_{2} & 0 & t_{1} \\
t_{1} & t_{2} & t_{0}
\end{pmatrix}
\nonumber \\
&=&
\begin{pmatrix}
t_{4} & 2t_{2} & t_{3} \\
t_{3} & t_{2} & t_{2} \\
t_{2} & 0 & t_{1}
\end{pmatrix}
=\begin{pmatrix}
2 & 2 & 1 \\
1 & 1 & 1 \\
1 & 0 & 0
\end{pmatrix},
\eea

and likewise
\begin{equation}
Q^n_{3} = 
\begin{pmatrix}
t_{n+2}[S_2(3)] & t_{n}[S_X(3)]~t_2[S_2(3)] & t_{n+1}[S_2(3)] \\
t_{n+1}[S_2(3)] & t_{n-1}[S_X(3)]~t_2[S_2(3)] & t_{n}[S_2(3)] \\
t_{n}[S_2(3)] & t_{n-2}[S_X(3)]~t_2[S_2(3)] & t_{n-1}[S_2(3)]
\end{pmatrix}
~~(n \ge 2),
\end{equation}
where the notations have their usual meanings (see Table \ref{table-p3}).\\

Hence, the determinantal identity is
\bea 
(\text{det}~Q_3)^n = 
\begin{vmatrix}
t_{n+2}[S_2(3)] & t_{n}[S_X(3)]~t_2[S_2(3)] & t_{n+1}[S_2(3)] \\
t_{n+1}[S_2(3)] & t_{n-1}[S_X(3)]~t_2[S_2(3)] & t_{n}[S_2(3)] \\
t_{n}[S_2(3)] & t_{n-2}[S_X(3)]~t_2[S_2(3)] & t_{n-1}[S_2(3)]
\end{vmatrix}
= 1. \label{eq-3cassini}
\eea

\subsection{Illustration for $p=4$}
For $p=4$, we have
\begin{equation}
Q_{4} = 
\begin{pmatrix}
1 & 1 & 1 & 1 \\
1 & 0 & 0 & 0\\
0 & 1 & 0 & 0 \\
0 & 0 & 1 & 0 \\
\end{pmatrix}
=
\begin{pmatrix}
t_{4} & t_{3} & t_{3} & t_{3} \\
t_{3} & 0 & 0 & t_{2} \\
t_{2} & t_{3} & 0 & t_{1} \\
t_{1} &0 &  t_{3} & t_{0}
\end{pmatrix},
\end{equation}
\bea
Q^2_{4} &=& 
\begin{pmatrix}
1 & 1 & 1 & 1 \\
1 & 0 & 0 & 0\\
0 & 1 & 0 & 0 \\
0 & 0 & 1 & 0 \\
\end{pmatrix}
\begin{pmatrix}
t_{4} & t_{3} & t_{3} & t_{3} \\
t_{3} & 0 & 0 & t_{2} \\
t_{2} & t_{3} & 0 & t_{1} \\
t_{1} &0 &  t_{3} & t_{0}
\end{pmatrix}
\nonumber \\
&=&
\begin{pmatrix}
t_{5} & 2t_{3} & 2t_{3} & t_{4} \\
t_{4} & t_{3} & t_{3} & t_{3} \\
t_{3} & 0 & 0 & t_{2} \\
t_{2} & t_{3} & 0 & t_{1}
\end{pmatrix}
=
\begin{pmatrix}
2 & 2 & 2 & 1 \\
1 & 1 & 1 & 1 \\
1 & 0 & 0 & 0 \\
0 & 1 & 0 & 0
\end{pmatrix},
\eea
and, for $n \ge 3$,
\begin{equation}
Q^n_{4} = 
\begin{pmatrix}
t_{n+3}[S_3(4)] & t_{n+3}[S_2(4)]~t_{3}[S_3(4)] & t_{n}[S_X(4)]~t_{3}[S_3(4)] & t_{n+2}[S_3(4)] \\
t_{n+2}[S_3(4)] & t_{n+2}[S_2(4)]~t_{3}[S_3(4)] & t_{n-1}[S_X(4)]~t_{3}[S_3(4)] & t_{n+1}[S_3(4)] \\
t_{n+1}[S_3(4)] & t_{n+1}[S_2(4)]~t_{3}[S_3(4)] & t_{n-2}[S_X(4)]~t_{3}[S_3(4)] & t_{n}[S_3(4)] \\
t_{n}[S_3(4)]     & t_{n}[S_2(4)]~t_{3}[S_3(4)] & t_{n-3}[S_X(4)]~t_{3}[S_3(4)] & t_{n_1}[S_3(4)]
\end{pmatrix},
\end{equation}
where the notations have their usual meanings (see Table \ref{table-p4}).\\

Hence, the determinantal identity for \(4\)-sequence $S_3(4)$ is
\bea 
(\text{det}~Q_4)^n &=& 
\begin{vmatrix}
t_{n+3}[S_3(4)] & t_{n+3}[S_2(4)]~t_{3}[S_3(4)] & t_{n}[S_X(4)]~t_{3}[S_3(4)] & t_{n+2}[S_3(4)] \\
t_{n+2}[S_3(4)] & t_{n+2}[S_2(4)]~t_{3}[S_3(4)] & t_{n-1}[S_X(4)]~t_{3}[S_3(4)] & t_{n+1}[S_3(4)] \\
t_{n+1}[S_3(4)] & t_{n+1}[S_2(4)]~t_{3}[S_3(4)] & t_{n-2}[S_X(4)]~t_{3}[S_3(4)] & t_{n}[S_3(4)] \\
t_{n}[S_3(4)]     & t_{n}[S_2(4)]~t_{3}[S_3(4)] & t_{n-3}[S_X(4)]~t_{3}[S_3(4)] & t_{n_1}[S_3(4)]
\end{vmatrix}
\nonumber \\~\nonumber \\
&=& (-1)^n. \label{eq-4cassini}
\eea

\end{document}